\documentclass[12pt,a4paper]{amsart}
\usepackage{amssymb,color,currfile,pdfpages}

\usepackage[english]{babel}
\usepackage{verbatim,here}
\usepackage[T1]{fontenc}
\usepackage{epstopdf}
\usepackage{floatflt,graphicx}
\usepackage{color,xcolor}
\usepackage{orcidlink}
\usepackage{hyperref}
\usepackage{subcaption}
\usepackage{url}

\usepackage{enumerate}
\usepackage{animate}
\usepackage{a4wide}
\usepackage{amsmath,amssymb}
\usepackage{amsthm}
\usepackage{float}
\usepackage{booktabs}
\usepackage{algorithm}
\usepackage{algorithmic}

\usepackage{tikz}
\usetikzlibrary{calc, math, angles}
\usetikzlibrary{decorations.pathreplacing,calligraphy, positioning}

\theoremstyle{remark}

\newcounter{minutes}
\divide\time by 60
\newcounter{hours}
\multiply\time by 60 \addtocounter{minutes}{-\time}
\newcommand{\bea}{\addtocounter{num}{1}\begin{eqnarray}}  
\newcommand{\eea}{\end{eqnarray}}
\dedicatory{}
\commby{}
\theoremstyle{plain}
\newtheorem{thm}[equation]{Theorem}

\theoremstyle{definition}

\theoremstyle{remark}

\newtheorem{nonsec}[equation]{}
\numberwithin{equation}{section}

\newcommand{\beq}{\begin{equation}}
\newcommand{\eeq}{\end{equation}}
\newcommand{\ben}{\begin{enumerate}}
\newcommand{\een}{\end{enumerate}}
\newcommand{\bequu}{\begin{eqnarray*}}
\newcommand{\eequu}{\end{eqnarray*}}
\newcommand{\bequ}{\begin{eqnarray}}
\newcommand{\eequ}{\end{eqnarray}}

\newcommand{\UH}{\mathbb{H}^2}
\newcommand{\BB}{\mathbb{B}^2}

\renewcommand{\Im}{{ \rm Im}\,}
\renewcommand{\Re}{{ \rm Re}\,}

\DeclareMathOperator{\dist}{dist}

\begin{document}
\thispagestyle{empty}
\def\thefootnote{}

%\title[Shortest path in planar domains]{Shortest paths in planar domains with hyperbolic type metrics}
%\title[Approximating geodesics of hyperbolic type metrics]{Approximating geodesics of hyperbolic type metrics on planar domains}

\title[Shortest paths of hyperbolic type metrics]{Shortest paths in planar domains with hyperbolic type metrics}

\author[S. Gao]{Shuliang Gao \orcidlink{0009-0000-5819-1717}}
\address{ Department of Mathematics with Computer Science, Guangdong Technion - Israel Institute of Technology, 241 Daxue Road, Jinping District, Shantou, Guangdong 515063, People's Republic of China and Department of Mathematics, Technion - Israel Institute of Technology, Haifa 32000, Israel
\newline
\href{https://orcid.org/0009-0000-5819-1717}{{\tt https://orcid.org/0009-0000-5819-1717}}
}
\email{gao09228@gtiit.edu.cn}

\author[A. Hakanen]{Anni Hakanen \orcidlink{0000-0001-7473-2456}}
\address{Turku Collegium for Science, Medicine and Technology (TCSMT), University of Turku, Finland and Department of Mathematics and Statistics, University of Turku, Finland\newline
\href{https://orcid.org/0000-0002-1734-8228}{{\tt https://orcid.org/0000-0001-7473-2456}}
}
\email{anehak@utu.fi}

\author[A. Rasila]{Antti Rasila$^{*}$  \orcidlink{0000-0003-3797-942X}}
\address{
Department of Mathematics with Computer Science, Guangdong Technion - Israel Institute of Technology, 241 Daxue Road, Jinping District, Shantou, Guangdong 515063, People's Republic of China and Department of Mathematics, Technion - Israel Institute of Technology, Haifa 32000, Israel\newline
\href{https://orcid.org/0000-0003-3797-942X}{{\tt https://orcid.org/0000-0003-3797-942X}}
}
\email{antti.rasila@iki.fi; antti.rasila@gtiit.edu.cn}

\author[M. Vuorinen]{Matti Vuorinen \orcidlink{0000-0002-1734-8228}}
\address{Department of Mathematics and Statistics,
University of Turku, Finland\newline
\href{https://orcid.org/0000-0002-1734-8228}{{\tt https://orcid.org/0000-0002-1734-8228}}
}
\email{vuorinen@utu.fi}

\date{\today}

%%%%%%%%%%%%%%%%%%%%%%%%%%%%%%%%%%%%%
\begin{abstract}
%We study the problem of finding the shortest path between given two points
%in a planar domain....
%The problem of finding a shortest path is very common in many situations in everyday life as well as in mathematical research. 
We study planar domains $G$ equipped with a hyperbolic type metric and approximate geodesics that join two points $x,y \in G$ and their lengths. We present an algorithm that enables one to approximate the shortest distance in polygonal domains taken with respect to the quasihyperbolic metric. The method is based on Dijkstra's algorithm, and we give several examples demonstrating how the algorithm works and analyze its accuracy. We experimentally demonstrate several previously theoretically observed features of geodesics, such as the relationship between hyperbolic and quasihyperbolic distance in the unit disk. We also investigate bifurcation of geodesics and the connection of this phenomenon to the medial axis of the domain.
\end{abstract}

%%%%%%%%%%%%%%%%%%%%%%%%%%%%%%%%%%%%%

\maketitle

\phantomsection
\footnotetext{${}^*$ Corresponding author. 

Matlab source code for this paper is available from:\\ \href{https://github.com/arasila/quasihyperbolic_geodesics}{{\tt https://github.com/arasila/quasihyperbolic\_geodesics}}\\
\texttt{{\tiny File:~\jobname .tex, printed: \number\year-%
\number\month-\number\day, \thehours.\ifnum\theminutes<10{0}\fi\theminutes}}}
\makeatletter

\makeatother

% Text of article.

%%%%%%%%%%%%%%%%%%%%%%%%%%%%%%%%%%%%%
%%%%%%%%%%%%%%%%%%%%%%%%%%%%%%%%%%%%%
%%%%%%%%%%%%%%%%%%%%%%%%%%%%%%%%%%%%%
\section{Introduction}

The problem of finding a shortest path between two given locations is a classical problem in geodesy (see e.g. \cite[Section 13.6]{VermeerRasila00}), and it has many connections to real world questions such as route planning, logistics, robotics, aviation, engineering, and communication networks. In this paper, we investigate this problem for hyperbolic type weighted metrics in planar domains (see \cite{hkv}). One of the methods to study these problems is to discretize the underlying mathematical model by using graphs. The vertices of the graph are possible intermediate locations between initial and terminal locations, whereas the edges depict all the connections between vertices. The research literature of this part of mathematical graph theory is vast due to the possible savings as a result of process speedup of algorithms.

We study this shortest path problem in a planar domain with a given metric with some
special features: in particular, this metric measures not only the distance between two points
but also their location with respect to the boundary of the domain. In the cases we study, the
existence of shortest paths between two given points, i.e. length-minimizing curves or geodesics, is known by theory. A characteristic feature of the geodesics of hyperbolic type metrics is that they have a tendency to avoid the boundary of the domain as much as possible. 
%Metrics of this type are sometimes called hyperbolic type metrics \cite{hkv}.

\subsection{Weighted metrics}
For a domain $G \subsetneq \mathbb{R}^n, n\ge 2,$ define a continuous function $w: G \to (0,\infty),$ called {\it weight function} and the {\it weighted
length} of a rectifiable curve $\gamma$ in $ G$ by
\begin{equation}
\label{w-len}
\ell_G^w(\gamma) = \int_{\gamma} \,w(z)\, |dz|\,.
\end{equation}
If $w(x)=1$ for all $x \in G\,,$ then $\ell_G^w(\gamma)$ reduces to the Euclidean length of
$\gamma$, and if 
\begin{equation} \label{qh-len}
w(z)= \frac{1}{d(z, \partial G)}\,,
\end{equation}
where $d(z, \partial G)$ is the distance of $z\in G$ from the boundary of the domain, then $\ell_G^w(\gamma)$ is the {\it quasihyperbolic
length} of $\gamma\,.$
The quasihyperbolic distance 
 is defined for $x,y \in G$ as
\begin{equation}\label{qhdef}
k_G(x,y)= \inf_{\gamma} \int_{\gamma} \, \frac{|dz|}{d(z, \partial G)}\,.
\end{equation}
where the infimum is taken over all rectifiable curves $\gamma$ in $G$ with
$x,y \in \gamma.$  This metric was defined by F.W. Gehring and his coauthors 
B.P. Palka \cite{gp} and B. Osgood \cite{go}, who also studied its applications in geometric function theory. In particular, they showed that the distance function \eqref{qhdef} defines a metric and that there exists a length
minimizing curve, geodesic, joining the points $x$ and $y\,.$ Interestingly, this is not necessarily true in more general settings like Banach spaces (see \cite{RasilaTalponen14}). 

%Moreover, the quasihyperbolic metric shares several properties of the hyperbolic metric and therefore is an example of a hyperbolic type metric.

%
%It is well known that this infimum is attained, i.e. the infimum is a minimum~\cite{go,gp}.
Since its introduction, the quasihyperbolic metric has become a standard tool in the study 
of geometric function theory,  quasiconformal mappings,  and  analysis on metric spaces \cite{gh, gmp, hkv, hei, vjussi}, where it is typically used as the substitute for hyperbolic metric in the cases where hyperbolic metrics cannot be defined. For half-planes and half-spaces, the quasihyperbolic metric coincides with the hyperbolic metric, and when the domain $G$ is a simply connected domain of the plane, it differs from the hyperbolic metric \cite{b} at most by a constant factor \cite{gh}. This serves as a motivation for the Gromov hyperbolic uniformization theory, introduced by Bonk, Heinonen, and Koskela in the seminal paper \cite{BonkHeinonenKoskela01}, where the question of Gromov hyperbolicity of quasihyperbolic metrics plays a central role. As far as we know, the algorithm given in this paper is the first general method for approximation of geodesics and distances in quasihyperbolic metric.

\begin{nonsec}{\bf Approximation of quasihyperbolic geodesics.}
In spite of all these applications, finding the explicit value of the quasihyperbolic distance is not possible
except in a few special cases.  Therefore, during the past five decades many comparison metrics  have been introduced  that are easily evaluated and, at the same time,  give either upper or lower bounds for the quasihyperbolic metric.  For a survey of these comparison metrics, see \cite{hkv} and the recent  article \cite{mo}. 
\end{nonsec}

In this study, quasihyperbolic geometry is investigated from the point of view of constructive approximation.
We discretize the problem by forming a grid of points in the domain $G\,.$
This grid must be fine enough considering the mutual locations of the points
and the geometric features of the domain $G\,.$  We define a graph using
these grid points and define the distance between two neighbouring points
in a manner that is compatible with the quasihyperbolic geometry.  The main
step of the process is to find the shortest path in this graph.  For this purpose 
we use Dijkstra's algorithm \cite{Dijkstra59}. This method applies not only for the
quasihyperbolic metric but also for many other weighted metrics for which shortest
paths joining the points are known to exist. The method enables us to approximate the
quasihyperbolic geodesics by means of polygonal curves and present figures of these approximations in several polygonal
domains. We also compare our approximate values
of the quasihyperbolic distances to approximations
obtained by other methods.

After this introduction the contents of this paper is organized as follows. In Section 2 we introduce some basic notation for the metrics we use, describe how to  proceed from grid points of a domain to
a graph, and some facts about special functions. In Section 3 we discuss computing distances in graphs and introduce Dijkstra's algorithm and discuss implementation of the algorithm for approximation of geodesics  based on Dijkstra's algorithm. In Section 4 we use this algorithm for approximation of hyperbolic and quasihyperbolic geodesics in the unit disk, with an aim of validating the method and understanding its accuracy. We also consider conformal mappings of the upper half-plane onto two polygonal domains, expressed in terms of Schwarz-Christoffel integrals. In Section 5 we investigate cases where non-uniqueness and bifurcation of geodesics is expected from theoretical considerations. In Sections 6--8 we consider other interesting examples where hyperbolic and quasihyperbolic geodesics can be approximated by using our methods. Finally, in Section 9, we give our final conclusions.
%%%%%%%%%%%%%%%%%%%%%%%%%%%%%%
%%%%%%%%%%%%%%%%%%%%%%%%%%%%%%
%%%%%%%%%%%%%%%%%%%%%%%%%%%%%%
\section{Preliminaries}
%%%%%%%%%%%%%%%%%%%%%%%%%%%%%%
%%%%%%%%%%%%%%%%%%%%%%%%%%%%%%
%%%%%%%%%%%%%%%%%%%%%%%%%%%%%%
In this section, we give the formulas and results needed later in our examples. We first focus on hyperbolic geometry and then present the necessary tools in order to handle conformal mappings of the upper half-plane onto a quadrilateral domain.

\begin{nonsec}{\bf Hyperbolic geometry.}\label{hg}
We recall some basic formulas and notation for hyperbolic geometry from \cite{b}.
The hyperbolic metrics of the unit disk ${\mathbb{B}^2}$ and
the upper half-plane  ${\mathbb{H}^2}$ are defined, respectively, by
\begin{equation}\label{rhoB}
\sinh \left( \frac{\rho_{\mathbb{B}^2}(a,b)}{2} \right)=
\frac{|a-b|}{\sqrt{(1-|a|^2)(1-|b|^2)}} ,\quad a,b\in \mathbb{B}^2\,,
\end{equation}
and
\begin{equation}\label{rhoH}
\cosh \left( \rho_{\mathbb{H}^2}(x,y)\right) = 1+\frac{|x-y|^2}{2 {\rm Im}(x) {\rm Im}(y)}
\, ,\quad x,y\in \mathbb{H}^2\,.
\end{equation}
By the Riemann mapping theorem, one can map a simply connected planar 
domain $D$, not equal to
the whole plane, onto the unit disk by means of a conformal mapping $f:D \to \mathbb{B}^2$  and define the hyperbolic metric for
 $a,b \in D$
as follows
\[
\rho_D(a,b)= \rho_{\mathbb{B}^2}(f(a),f(b))  \,. 
\]
This definition is possible due to the \emph{conformal invariance} of the hyperbolic
metric \cite{kr}.

Numerical methods to find the hyperbolic distance between two given
points in a simply connected planar domain and to draw the hyperbolic geodesic joining the points are given in \cite{NV}.

\end{nonsec}
%\includepdf[pages=-,pagecommand={},width=\textwidth]%{excentricity20230427.pdf}

\begin{nonsec}{\bf Distance ratio metric.}
 A useful comparison
metric for the quasihyperbolic metric is the distance ratio metric \cite{gp} 
\begin{equation}\label{jdef}
j_G(x,y)= \log \left( 1+ \frac{|x-y|}{\min\{ d(x, \partial G),  d(y, \partial G)\}}\right) \,,
\end{equation}
where $x,y \in G$.
The metrics $j_G$ and $k_G$ are not conformally invariant, and they are only 
similarity invariant.
If $G \in  \{  \mathbb{B}^2,  \mathbb{H}^2\},  $ then for
$x,y  \in G$
 \begin{equation} \label{jrho}
 j_G(x,y)  \le \rho_G(x,y) \le 2 j_G(x,y) \,,
      \end{equation}
 \begin{equation} \label{krho}
 k_G(x,y)  \le \rho_G(x,y) \le 2 k_G(x,y) \,.
      \end{equation}      
In fact, $k_G \equiv \rho_G $ for $G= \mathbb{H}^2.$
The geodesics of the quasihyperbolic metric are smooth
\cite{mar}.

For all domains $G,$ we have a simple lower bound for the quasihyperbolic metric
\begin{equation} \label{kjineq}
k_G(x,y)\ge  j_G(x,y) \ge \left|  \log  \frac{d(x, \partial G)}{  d(y, \partial G)} \right|\,.
\end{equation}
\end{nonsec}

  In view of the above two-sided inequalities  \eqref{jrho} and  \eqref{krho}
and the lower bound  \eqref{kjineq}  it is natural to ask for which classes of
domains there is a majorant for $k_G$ in terms of $j_G.$
  Indeed, inequalities between
metrics can  be used to define some important classes of domains.  Thus,  for instance, the class of simply connected plane domains $G$ for which there exists a constant $A \ge 1$ such that for all $x,y \in G$
\begin{equation} \label{unifdef}
k_G(x,y) \le A \, j_G(x,y)
\end{equation}
defines the class of quasidisks studied in \cite{gh}.  Domains of this class, in turn, are
special cases of the so called $\varphi$-uniform domains \cite{hkv}.

\begin{nonsec}{\bf The medial axis.}
The medial axis of a planar domain is the set of all those points of the domain that have at least two closest boundary
points. Thus, for a sector with angle $\theta \in (0, \pi),$ the angle bisector is the medial axis. For an ellipse domain,
a part of the major axis is the medial axis.

Lind\'en showed in his PhD thesis \cite{Linden} that if two points of a convex sector  are far enough from each other,
then the quasihyperbolic geodesic joining them consists of  three pieces linked together.
Two pieces are sub-arcs of circles, each circle orthogonal to a side of the sector and tangential to the sector bisector at 
a contact point, whereas the third piece is the subsegment  of the bisector  joining these two contact points. See also H\"ast\"o \cite{hasto}.

\end{nonsec}

\begin{nonsec}{\bf Approximation of the weighted length of a polygonal curve. } Let $x_j, j=1,\ldots,p,$ be points in a domain
$G \subset \mathbb{R}^n$ and consider the polygonal curve
$\gamma= \cup_{j=1}^{p-1} \, [x_j, x_{j+1}]\,.$ For a given weight function $w:G \to (0,\infty)$ we can approximate
$\ell_G^w (\gamma)$ provided that $\gamma \subset G$ which clearly holds if $G$ is convex or , more generally, if
all the segments $ [x_j, x_{j+1}] \subset G, j=1,\ldots,p-1\,.$
In what follows we assume that this last condition is fulfilled.

From the definition \eqref{qhdef} of the quasihyperbolic metric
it is clear that a simple approximation is
\begin{equation}
    \label{lenapp1}
    \ell_G^w (\gamma) \approx \sum_{j=1}^{p-1}\, w_j|x_j - x_{j+1}|
\end{equation}
where the coefficient $w_j$ depends on the weight function
$w\,.$ It is a natural requirement that $w_j$ depends on the
points $x_j$ and $x_{j+1}$ in a symmetric manner, because $\ell_G^w(x_j, x_{j+1}) = \ell_G^w(x_{j+1}, x_j) \,.$
Therefore we can choose e.g.
\begin{equation}
    \label{qhwt}
    w_j = \frac{1}{\min \{ d(x_j, \partial G),d(x_{j+1}, \partial G) \}}
\end{equation}
in the case of the quasihyperbolic metric; this is what we use in the experiments reported below.

The case of a general weight function is more difficult because the oscillation of the function $w$ is uncontrolled.
For the purpose of controlling the oscillation of the function
$w$ one could require that $w$ satisfies a Harnack inequality:
there exist constants $\lambda\in (0,1)$ and $C_{\lambda}\ge 1$ such that
\begin{equation}
    \label{harnineq}
    \max_{B_x} w(z) \le C_{\lambda} \,  \min_{B_x} w(z)
\end{equation}
holds for all $B^n(x,r)\subset G $ where $B_x= \overline{B}^n(x, \lambda r)\,.$ If it is true that, for the above
polygon and for all $j=1,\ldots,p-1,$ 
$$[x_j, x_{j+1}] \subset B^n(x_j, \lambda d(x_j))$$
and the Harnack condition holds, then we could estimate
the error in the approximation \eqref{lenapp1}. Note that the
quasihyperbolic weight function $1/d(x,\partial G) $ satisfies
the above Harnack condition.
\end{nonsec} 

\begin{nonsec}{\bf Grid points.} \label{grpts} For the purpose of
approximating distances between point pairs in a
planar domain $G\,,$ we generate, for a given mesh
size $h>0,$ a grid of those points of the
Cartesian lattice $h \mathbb{Z}^2,$ which are in the domain $G\,.$

Then for each grid point we define its set of neighboring points. We use the geographical directions to label these, the four cardinal directions -- north (N), east (E), south (S), and west (W) -- and the four intermediate directions -- northeast (NE), southeast (SE), southwest (SW), and northwest (NW), see Figure \ref{fig:mesh_directions}. The set of neighboring
points of a grid point $x$ may either be empty
or contain $p$ points, $1\le p\le 8 \,,$ depending
on $d(x, \partial D)$ and the mesh size $h\,.$ 
Observe that if $x$ is a neighbor of $y\,$ then
also $y$ is a neighbor of $x\,.$
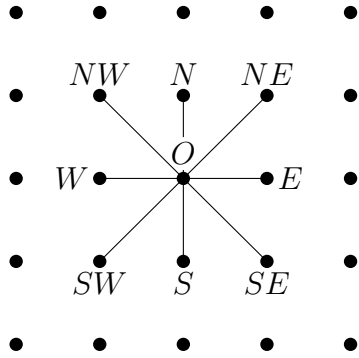
\begin{figure}
	\centering
	\begin{tikzpicture}[scale=1.1]
		%--BASIC SHAPE---------------------
		%------DEFINITIONS-----------------
		%-------LINES-----------------------
		\draw \foreach \x in {(1,0),(0,1),(-1,0),(0,-1),(1,1),(1,-1),(-1,1),(-1,-1)} {
			(0,0) -- \x
		};
		%-------NODES-----------------------
		\draw \foreach \x in {-2,-1,0,1,2} \foreach \y in {-2,-1,0,1,2} {
			(\x,\y) node[circle, fill=black, inner sep=0pt, minimum width=5pt] {} }; % black nodes 
		%--END OF BASIC SHAPE---------------
		%-----------------------------------
		%--LABELS--------------------
		\draw {
			(0,.3) node[circle,fill=white,inner sep=0pt] {$O$}
			(1,0) node[right] {$E$}
			(0,1) node[above] {$N$}
			(-1,0) node[left] {$W$}
			(0,-1) node[below] {$S$}
			(1,1) node[above] {$NE$}
			(1,-1) node[below] {$SE$}
			(-1,1) node[above] {$NW$}
			(-1,-1) node[below] {$SW$}
		};
	\end{tikzpicture}
\caption{A grid point O with a maximal number, eight,  neighboring points as described in \ref{grpts} (without a refinement of the grid).}\label{fig:mesh_directions}
\end{figure}

\end{nonsec}

\begin{nonsec}{\bf From a grid to a graph. } 
We now set up a graph using the grid points of a domain
$G$ constructed above as vertices of a graph and defining the distance $d(x,y)$ between two
neighboring points $x$ and $y$ such that $d(x,y) = d(y,x).$ In what follows, $d$ is usually a
weighted metric, in most cases the quasihyperbolic metric.
For the purpose of approximating the shortest path between two grid points $x$ and $y,$
we will apply Dijkstra's algorithm. For our method to work, there must exist at
least one polygonal curve $\cup^{k-1}_{j=1} 
 [p_j ,p_{j+1}  ],$ with $p_1
= x$ and $p_k = y$ such that
$p_{j+1}$  is a neighbor of $p_j$ and $[p_j ,p_{j+1}  ] \subset G$, for each $j = 1,\ldots,k-1.$ It is easily seen that
this last condition is satisfied if $d(p_j ,\partial G) \ge  2h,$ for all $j = 1,\ldots,k .$ Dijkstra's algorithm
provides a simple and robust method for finding the shortest distance between two
points of the graph.
\end{nonsec}

\begin{nonsec}{\bf From a graph to an algorithm.}   The above method can be improved by refining the
grid. We now describe the details of this  refinement idea. 
\begin{comment} OLD VERSION:
    Based on Subsection \ref{grpts} we improve the concept of a neighborhood with the idea of concentric circles. The set of neighboring points of a grid point $x$ is defined with a positive parameter $m \in \mathbb{N}$ as follows:
\begin{equation}\label{eq: ngbh_set}
    N_x=\big\{p \in G\;\big|\;h\cdot m\leq |p-x| \leq \sqrt 2 \cdot h\cdot m\big\}.
\end{equation}
\end{comment}
Indeed, we increase the number of neighboring points by accepting
certain points of the refined grid $G \cap \frac{h}{m} \mathbb{Z}^2$ where  $m \in \mathbb{N},$ as neighbors of a given point $x$ in the original grid. This expanded set of neighbors is defined  for
 a positive parameter $m \in \mathbb{N}$ as follows:
\begin{equation}\label{eq: ngbh_set}
    N_x=\{p \in G  \cap \frac{h}{m} \mathbb{Z}^2\;|\;h\cdot m\leq |p-x| \leq \sqrt 2 \cdot h\cdot m\}.
\end{equation}
Figure \ref{fig:ngbh_points} illustrates the set $N_x$ for different values of $m$. Observe that if $N_x$ contains $y$, then also $N_y$ contains $x$.
\begin{comment} OLD VERSION:
 Although $h$ represents the mesh size and $m$ influences the direction, we have not standardized $h$ and $m$ uniformly in this paper. This omission, however, does not imply that such standardization is infeasible. We plan to implement this operation in our future research.
 \end{comment} 
 This expansion of neighboring points means that 
the cardinality of the set $N_x$ depends on $m$ and can, in particular, be larger than $8.$ We have not fixed the parameter
$m$ uniformly in this paper and plan to study in our future
research the influence of $m$ on the accuracy of the computation.

\begin{figure}[htbp]
    \centering
    
    \begin{subfigure}[b]{0.45\textwidth}
        \centering
        \includegraphics[width=\linewidth]{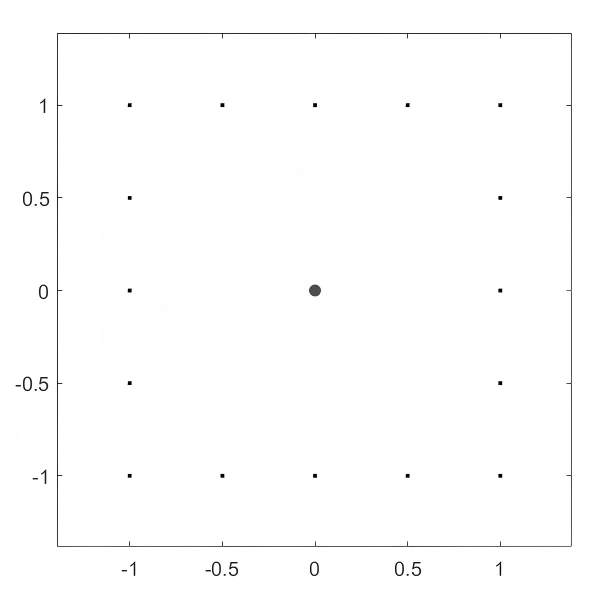}
        \caption{$m=2$}
        \label{fig:Layer_2}
    \end{subfigure}
    \hfill
    \begin{subfigure}[b]{0.45\textwidth}
        \centering
        \includegraphics[width=\linewidth]{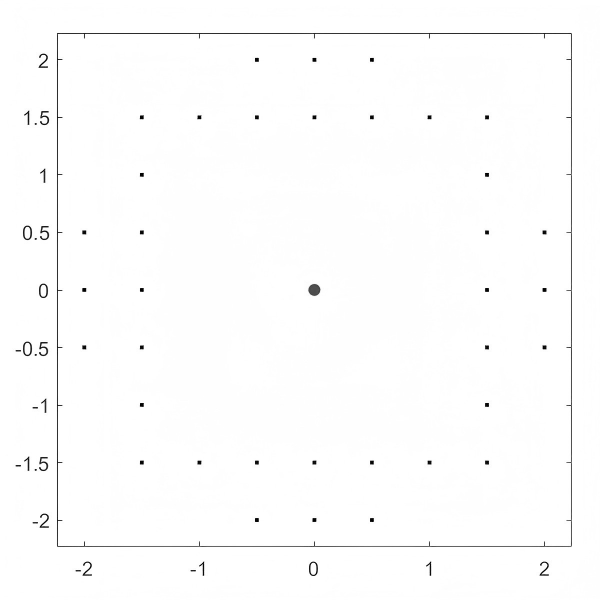}
        \caption{$m=3$}
        \label{fig:Layer_3}
    \end{subfigure}
    
    \begin{subfigure}[b]{0.45\textwidth}
        \centering
        \includegraphics[width=\linewidth]{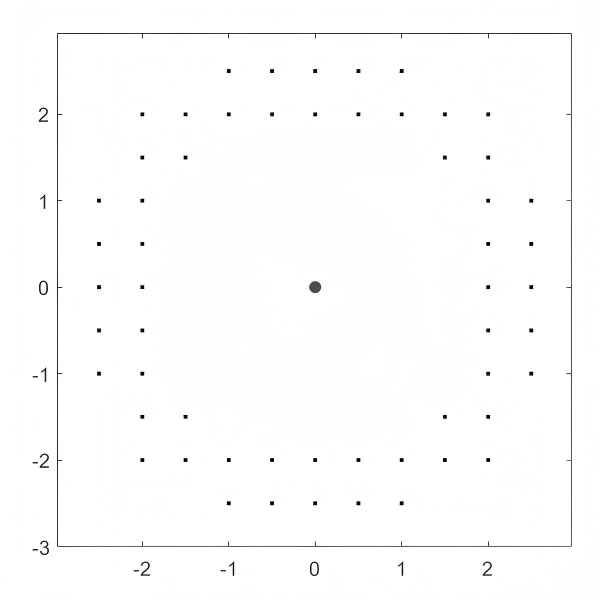}
        \caption{$m=4$}
        \label{fig:Layer_4}
    \end{subfigure}
    \hfill
    \begin{subfigure}[b]{0.45\textwidth}
        \centering
        \includegraphics[width=\linewidth]{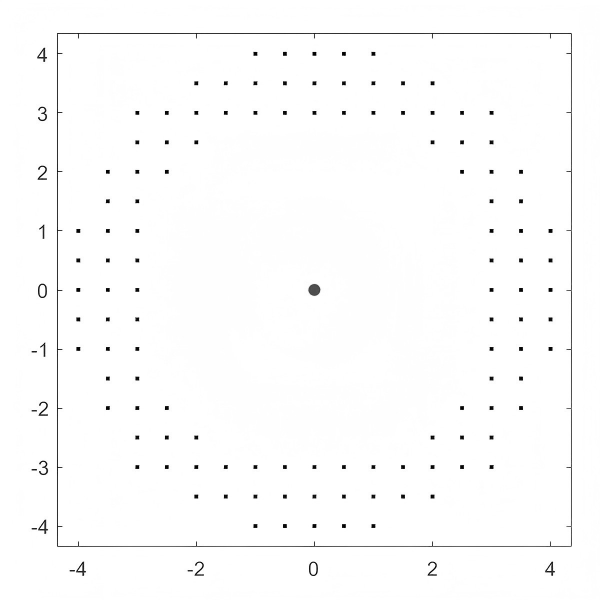}
        \caption{$m=6$}
        \label{fig:Layer_6}
    \end{subfigure}
    
    \begin{subfigure}[b]{0.45\textwidth}
        \centering
        \includegraphics[width=\linewidth]{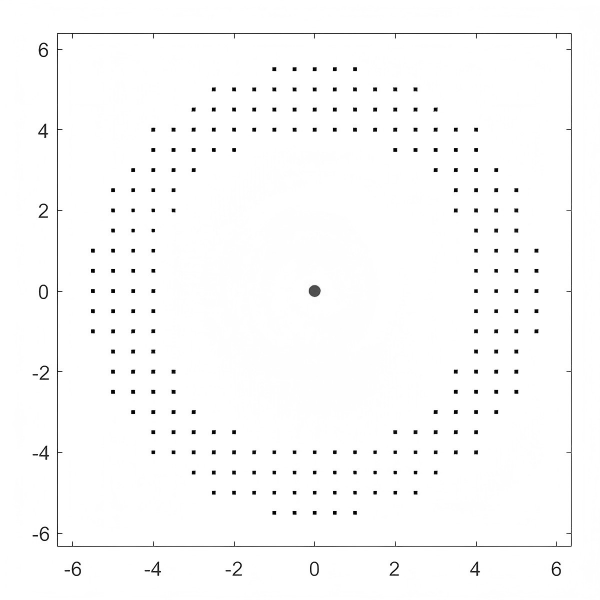}
        \caption{$m=8$}
        \label{fig:Layer_8}
    \end{subfigure}
    \hfill
    \begin{subfigure}[b]{0.45\textwidth}
        \centering
        \includegraphics[width=\linewidth]{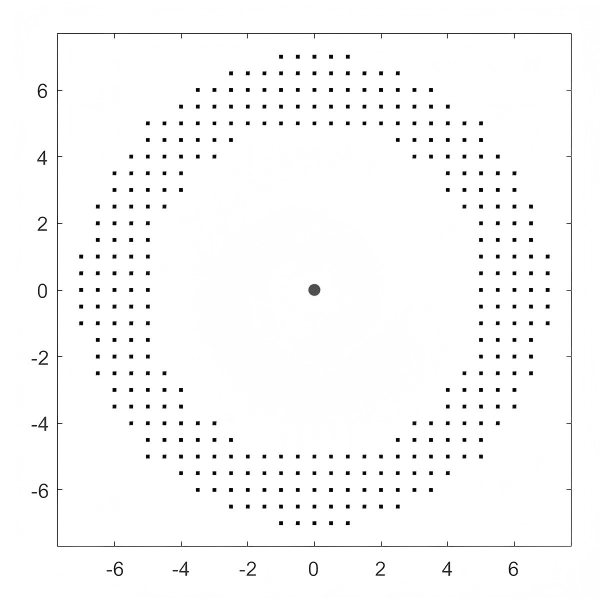}
        \caption{$m=10$}
        \label{fig:Layer_10}
    \end{subfigure}

    \caption{The $m$-neighborhoods of origin point $x=0$ with $h=0.5$.}
    \label{fig:ngbh_points}
\end{figure}

We now set up a graph using the grid points of a domain $G$ constructed above as vertices of a graph by using the method mentioned above to construct $N_x$ and to approximate the distance $d(x,y)$ between two neighboring points $x$ and $y \in N_x$ such that $d(x,y)= d(y,x)$. In what follows, $d$ is a weighted metric, in most cases the quasihyperbolic metric. 

We will apply Dijkstra's algorithm to approximate the shortest path between two grid points $x$ and $y$. For this method to work, there must exist at least one polygonal curve $\cup_{j=1}^{k-1}[p_j,p_{j+1}]\,,$
$p_j, j=1,\ldots,k,$ with $p_1=x$ and $p_k= y$
such that $p_{j+1}$ {\rm
is a neighbor of} $p_j$ 
and $[p_j, p_{j+1}] \subset G\,,$
for
each $j=1,\ldots,k-1\,.$ It is easily seen that
this last condition is satisfied if $d(p_j, \partial G) \ge 2 h\,,$
for
all $j=1,\ldots,k$. Our algorithm is described in detail in Section \ref{DM_Method}. 
\end{nonsec}

%\begin{nonsec}{\bf Mesh refinement.} The above method for finding the shortest path, while simple and robust, can be made more accurate through mesh refinement. For a given grid point $p$, the mesh size parameter $h$ and the direction parameter $m$, we introduce a refined grid of mesh size as follows. By using this definition, we may keep each step length bounded and add more directions while decreasing $h$ and increasing $m$, see Figure \ref{fig:ngbh_points}.
%$h/2$ and include suitable points from the refined
%grid as the neighbors of $p\,.$ A neighborhood candidate $q$ from the refined grid must satisfy
%at least $[p,q] \subset D\,.$ This refinement process can be iteratively repeated. 
%\end{nonsec}

\begin{nonsec}{\bf Special functions.}
Conformal mappings expressed in terms of special functions are frequently used in study of conformal invariants. In this paper, we apply the Schwarz-Christoffel formula which gives a conformal mapping of the upper half-plane onto a polygon in terms of an integral. This formula is useful in the case where the polygon is a quadrilateral and the conformal mapping is given by a hypergeometric function.

Let $\Gamma$ denote Euler's {\em gamma function} $\Gamma(z)$ and
$B(a,b)$ the beta function
\[ B(a,b)= \frac{\Gamma(a)\Gamma(b)}{\Gamma(a+b)}.\]
Given complex numbers
$a,b,$ and $c$ with $c\neq 0,-1,-2, \dots $,
the {\em Gaussian hypergeometric function} is the analytic
continuation to the slit plane $\mathbb{C} \setminus [1,\infty)$ of
the series
\begin{equation}
\label{eq:hypdef}
F(a,b;c;z) = {}_2 F_1(a,b;c;z) = 
\sum_{n=0}^{\infty} \frac{(a,n)(b,n)}{(c,n)} \frac{z^n}{n!}\,,\:\: |z|<1 \,.
\end{equation}
Here $(a,0)=1$ for $a \neq 0$, and $(a,n)$
is the {\em shifted factorial function}
or the {\em Appell symbol}
$$
(a,n) = a(a+1)(a+2) \cdots (a+n-1)
$$
for $n \in \mathbb{N} \setminus \{0\}$, where 
$\mathbb{N} = \{ 0,1,2,\ldots\}$.
Denote $C = C(b,c) = 1/B(b,c-b)$.
We define the function $f$
given on the closed upper half-plane by
\begin{eqnarray}  
w = f(z) & = &  f_{a,b,c}(z)  = C \int_{0}^{z} g_r(t) dt  \\ 
 & = & C \int_{0}^{z} 
t^{b-1}(1-t)^{c-b-1}(1-r^2t)^{-a} dt \label{eq:ellsin} \\ 
 & = &  e^{i(a+b+1-c)\pi} C r^{-2a} 
\int_{0}^{z} t^{b-1}(t-1)^{c-b-1}(t-1/r^2)^{(1-a)-1} dt\,. 
\end{eqnarray}
Recall the Euler integral representation \cite[Theorem 2.2.1]{aar}  and
\cite[15.3.1]{as}
\begin{equation} \label{eq:Fint}
F(a,b;c;z) = 
\frac{\Gamma(c)}{\Gamma(b)\Gamma(c-b)} %WAS on 040704:\frac{1}{ B (b) B(c-b)}  
\int_{0}^{1} t^{b-1}(1-t)^{c-b-1}(1-tz)^{-a} dt
\end{equation}
$$
= C(b,c) e^{i(a+b+1-c)\pi} \int_0^1 t^{b-1}(t-1)^{c-b-1}(tz-1)^{-a} dt
$$
% 2006-07-03
for ${\rm Re}(c) > {\rm Re}(b) > 0$ and
$z \in \mathbb{C} \setminus \{u \in \mathbb{R} \: | \: u \ge 1 \}$. % 2006-07-06
\end{nonsec}

\begin{thm} \cite[Theorem 2.4]{hvv}
\label{ellsintheo}
Let $H$ denote the closed upper half-plane $\{z \in \mathbb{C} \: | \: 
\Im (z) \ge 0\}$ and let  
$0 < a, b < 1, \max \{a+b,1 \} \le c \le 1 + \min \{a,b \}$,
$r \in (0,1)$. Then the
function $f$ in (\ref{eq:ellsin}) is a homeomorphism of $H$ onto the
quadrilateral $Q$ with vertices 
$$
f(0) = 0\,,\:\: f(1) = F(a,b;c;r^2)\,,
$$
$$
f(1/r^2) = f(1) + 
\frac{B(c-b,1-a)}{B(b,c-b)} e^{(b+1-c) i \pi} {(1-r^2)}^{2(c-a-b)}  
F(c-a,c-b;c+1-a-b;{(1-r^2)}^2)\,,
$$
and
$$
f(\infty) = f(1/r^2) 
+ \frac{B(1-a,a+1-c)}{B(b,c-b)}  e^{i(a+b+1-c)\pi} r^{2(1-c)} {r'}^{2(c-a-b)} 
F(1-b,1-a;2-c;r^2)\,,
$$
and interior angles $b\pi , (c-b)\pi, (1-a)\pi$ and $(a+1-c)\pi$, 
respectively, at these vertices. It is conformal in the interior of $H$.
\end{thm}

%%%%%%%%%%%%%%%%%%%%%%%%%%%%
%%%%%%%%%%%%%%%%%%%%%%%%%%%%
%%%%%%%%%%%%%%%%%%%%%%%%%%%%
\section{Dijkstra's Algorithm and Approximation of Geodesics}
\label{DM_Method}
%%%%%%%%%%%%%%%%%%%%%%%%%%%%
%%%%%%%%%%%%%%%%%%%%%%%%%%%%
%%%%%%%%%%%%%%%%%%%%%%%%%%%%

\begin{nonsec}{\bf Graphs and paths.}
A \emph{digraph} consists of vertices that are connected to each other by arcs. Formally, the vertex set of a digraph $D$ is denoted by $V(D)$ and the arc set by $A(D)$. The arc set $A(D)$ consists of ordered pairs $(u,v)$ where $u,v \in V(D)$. An undirected graph is a special case of a digraph where $(u,v) \in A(D)$ if and only if $(v,u) \in A(D)$. Digraphs do not necessarily satisfy this bidirectionality of arcs. A \emph{weighted digraph} is a digraph where each arc $(u,v)$ is associated with a \emph{weight} $w((u,v)) \in \mathbb{R}$.

A \emph{path} from a vertex $u_0$ to another vertex $u_k$ is a sequence of vertices $u_0 u_1 \cdots u_k$ where no vertex repeats and $(u_i,u_{i+1}) \in A(D)$ for all $u_i$. The distance from $u$ to $v$ (denoted by $\dist (u,v)$) is the minimum number of arcs along a path from $u$ to $v$ when the digraph is unweighted. If the digraph is weighted, the distance from $u$ to $v$ is defined as
\[ \dist (u,v) = \min\left\{\sum_{i=0}^{k} w((u_i,u_{i+1})) \, \big|\, u_0, \ldots, u_k \text{ is a path from } u =u_0 \text{ to } v=u_k \right\}. \]
\end{nonsec}

\begin{nonsec}{\bf Dijkstra's algorithm.}
The distances in a weighted digraph with non-negative weights can be computed with Dijkstra's algorithm \cite{Dijkstra59} which was introduced in 1959. This algorithm computes the distances from a single vertex $u$ to all other vertices in $D$. At first, each vertex $v\in V(D)$ is given an estimate $d(v)$ of the distance from $u$ to $v$. Initially, $d(u) = 0$ and $d(v)=\infty$ for all $v \in V(D) \setminus \{u\}$. 
Along the execution of the algorithm, we mark vertices and initially all vertices are unmarked. The following step is repeated until all vertices are marked.

We choose the vertex $x$ to be an unmarked vertex with $d(x)$ as small as possible ($x=u$ in the first iteration).
For all $v$ such that $(x,v) \in A(D)$, we update $d(v)$ to be $d(x) + w((x,v))$ if this is smaller than the current value of $d(v)$. Finally, we mark $x$ and sort the unmarked vertices according to the current estimates.

When this step is iterated until all vertices are marked, the execution ends and we have $d(v) = \dist (u,v)$ for all $v \in V(D)$. The sorting of the vertices according to their distance estimates can be done using Fibonacci heaps \cite{FredmanT87}. Then the running time of the algorithm is $\mathcal{O}(m + n \log n)$, where $m = |A(D)|$ and $n = |V(D)|$.
\end{nonsec}

\begin{nonsec}{\bf Dijkstra's algorithm compared to other path finding methods.}
Despite its age and simplicity, Dijkstra's algorithm remains one of the fastest and most used algorithms of computing distances in weighted graphs. If negative edge weights are allowed (but there are no negative cycles), the distances can be computed using the Bellman-Ford-Moore algorithm \cite{Bellman58,Ford56,Moore59}. If the digraph is unweighted, it is faster to use a simple algorithm adapted from breadth-first-search. During the preparation of this article, Duan et al. \cite{Duan25} showed that there exists an algorithm running in $\mathcal{O} (m \log^{2/3} n)$ time for computing the distances from one vertex to all others in a weighted digraph with non-negative edge weights. This is the long-awaited improvement to Dijkstra's algorithm, which has turned out difficult to improve due to the sorting step. However, Duan et al. make several assumptions to simplify the presentation of their algorithm. Most notably, the graphs are assumed to be sparse and every vertex has constant in- and out-degree. The algorithm they present is generalizable to graphs without these constraints. Dijkstra's algorithm is readily available and well-optimized for different programming languages, which is why we chose to use Dijkstra's algorithm and not to implement Duan et al.'s algorithm for the purposes of our experiments. 
\end{nonsec}
%This paper is very famous these days in China, since they avoid "sorting step" in traditional Dijkstra Algorithm. I have read that paper. In the paper they improve the efficiency on \textbf{sparse graph}, which has a definition below:

%\textbf{Def. Sparse Graph}

%Let $G=(V,E)$ be a graph. If $|E| \ll |V|^2$, then $G$ is a sparse graph.

%Our graph $G^*$ has $|E| < 8\cdot dir \cdot |V| \ll |V|$ for $|V| \to \infty$, so it is reasonable to apply method in this paper. However, Matlab already has built-in functions to apply DA to find shortest-path, if we want to apply the latest method, it will take at least a month for me to embed the new function into our code, i.e. , it will run out of our scheduled time. I hereby wait your order.

%%%%%%%%%%%%%%%%%%%%%%%%%%%%
%%%%%%%%%%%%%%%%%%%%%%%%%%%%
%%%%%%%%%%%%%%%%%%%%%%%%%%%%

%%%%%%%%%%%%%%%%%%%%%%%%%%%%
%%%%%%%%%%%%%%%%%%%%%%%%%%%%
%%%%%%%%%%%%%%%%%%%%%%%%%%%%

\begin{nonsec}{\bf Implementation of an algorithm for approximation of geodesics.}
  Due to practical concerns with  the numerical approximation, we use a version of the conformal density function modified from the usual definitions of hyperbolic and quasihyperbolic densities. 
\end{nonsec}

\begin{nonsec}{\bf Dijkstra's algorithm and weighted metrics.}
    We first generate an undirected weighted graph $G^*=(V,E,W)$ with mesh size $h$, parameter $m$ and definition of neighborhood points defined in \ref{eq: ngbh_set}, and then apply Dijkstra algorithm on $G^*$. The precise definition of the applied method is given as in Algorithm \ref{alg:combined_algorithm}.

\end{nonsec}

\small

\begin{algorithm}[htbp]
\caption{Graph Construction and Shortest Path Approximation}
\label{alg:combined_algorithm}
\begin{algorithmic}[1]
\REQUIRE Region $G \subset [a,b] \times [c,d]$, mesh size $h$, layer parameter $m$, weight function $w$, start point $p_1$, end point $p_2$
\ENSURE Undirected weighted graph $G^* = (V, E, W)$, shortest path $\gamma^*$ between $p_1$ and $p_2$, shortest path's length $\ell$

\STATE \textbf{Part 1: Graph Construction}
\STATE \textbf{Step 1:} Construct vertices $V$.
\STATE $V \gets \emptyset$
\FOR{$j = 0$ \TO $\lfloor (b-a)/h \rfloor$}
    \FOR{$k = 0$ \TO $\lfloor (d-c)/h \rfloor$}
        \STATE $x \gets a + j \cdot h\in\mathbb{R}$
        \STATE $y \gets c + k \cdot h\in\mathbb{R}$
        \STATE $p \gets x+ y \cdot i\in\mathbb{C}$
        \IF{$p \in G$}
            \STATE $V \gets V \cup \{p\}$
        \ENDIF
    \ENDFOR
\ENDFOR

\STATE \textbf{Step 2:} Construct edges $E$.
\STATE $E \gets \emptyset$
\FOR{ $p \in V$}
    \STATE $\text{Layer}_m(p) \gets \{ q \in V : m\cdot h \leq \|p - q\| \leq \sqrt{2} \cdot m\cdot h \}$
    \FOR{ $p_i \in \text{Layer}_m(p)$}
        \IF{Segment $[p,p_i] $ does not intersect with $G$}
            \STATE $E \gets E \cup \{(p, p_i)\}$
        \ENDIF
    \ENDFOR
\ENDFOR

\STATE \textbf{Step 3:} Construct weights $W$:
\STATE $W \gets \emptyset$
\FOR{ $e \in E$}
    \STATE $W \gets W \cup \{w(e)\}$
\ENDFOR
\STATE $G^* = (V, E, W)$

\STATE \textbf{Part 2: Shortest Path Calculation}
\STATE \textbf{Step 1:} Find the closest points to $p_1,p_2$ in $V$.
\STATE $\text{idx1} \gets \arg\min |zv - p_1|$ 
\STATE $\text{idx2} \gets \arg\min |zv - p_2|$

\STATE \textbf{Step 2:} Find shortest path and its length.
\STATE $[\gamma^*, \ell] \gets \textbf{shortestpath}(G^*, \text{idx1}, \text{idx2})$ \COMMENT{Dijkstra's algorithm returning both path and length\cite{matlab_shortestpath}}

\RETURN $G^*$, $\gamma^*$, $\ell$
\end{algorithmic}
\end{algorithm}

\normalsize
%\textbf{RMK 1}

\begin{nonsec}{\bf Remark.}
    It should be noted that a straight-forward application of the above methods leads to low efficiency because in most cases it is enough to carry out the refinement in a subdomain instead of the whole domain. Therefore, it is better to generate the uniform mesh in the region around $a^*, b^*$, defining $x_{\min},x_{\max},y_{\min},y_{\max}$ as below:
     %\textbf {Definition for $k$ (need to be complete)}
    \[
    \left \{ \begin{array}{cccc}
         x_{\min} & := & \min\{\Re(a^*), \Re(b^*)\} - k_{x_{\min}}\cdot h  ,\\
         x_{\max} & := & \max\{\Re(a^*), \Re(b^*)\} + k_{x_{\max}}\cdot h  ,\\
         y_{\min} & := & \min\{\Im(a^*), \Im(b^*)\} - k_{y_{\min}}\cdot h  ,\\
         y_{\max} & := & \max\{\Im(a^*), \Im(b^*)\} + k_{y_{\max}}\cdot h,
    \end{array}  \right.
    \]
    where $k_{x_{\min}},k_{y_{\min}},k_{x_{\max}},k_{y_{\max}}$ are set manually. Then we apply Algorithm \ref{alg:combined_algorithm} on the region $[x_{\min}, x_{\max}] \times [y_{\min}, y_{\max}] \cap G$ to generate the local approximated geodesic. We mainly apply this method in our following experiments, and call it applying Algorithm \ref{alg:combined_algorithm} locally.
\end{nonsec}

%\begin{rem}
    %Add the link that shows how does shortestpaths and distances are defined in MATLAB. 
    %\textbf{Here I've still not decided where to put this cite. Moreover, please check the format of the cite.}
%\end{rem}

%\begin{nonsec}{\bf Remark.}
    %From this section we use DM refer to Dijkstra Algorithm method. It means that when stating "apply DM on $G$", it is equivalent to apply \ref{alg:weighted_mesh_generation} on $G$ first, and then apply \ref{alg:shortest_path_calculation} to obtain the approximated geodesic $\gamma^*$.
%    \textbf{Replace all DM with "Algorithm \ref{alg:combined_algorithm} with refinement"}
%    From this section we use DM refer to Algorithm  \ref{alg:combined_algorithm}. It means that when stating "apply DM on $G$", it is equivalent to apply \ref{alg:combined_algorithm}. Moreover, when stating "doing refinement on DM", it is equivalent to decrease $h$ and increase $m$ reasonably.
%\end{nonsec}

\begin{nonsec}{\bf Hyperbolic geodesics.}
%The following are the constructions for different curves.
Let $G \subset \mathbb{C}$ be a simply connected domain, and let $f: \mathbb{H}^2 \to G$ be a conformal mapping (implemented via Schwarz-Christoffel transformation). For any $a^*, b^* \in G$, the discretization of the hyperbolic geodesic $\gamma_{\rho_G}(a^*,b^*)$ on $G$ is constructed as follows:
 
Let $a = f^{-1}(a^*)$, $b = f^{-1}(b^*)$. Then we know the hyperbolic geodesic $\gamma_{\rho_{\mathbb{H}^2}}(a,b)$ on $\mathbb{H}^2$. We choose a discrete set of points $\{z_i\}_{i=1}^n \subset \gamma_{\rho_{\mathbb{H}^2}}(a,b)$ satisfying:
    \begin{align*}
        z_1 &= a, \quad z_n = b, \\
        \|z_{i+1} - z_i\| &= \Delta s, \quad i = 1,\dots,n-1,
    \end{align*}
    where $\Delta s$ is the size of the discretization step. Then we obtain the point sequence on $G$ via conformal mapping: $w_i = f(z_i)$, $i = 1,\dots,n$. The final curve is defined as:
    %\[
    %\gamma_{\rho_G}(a^*,b^*) = \bigcup_{i=1}^{n-1} \overline{w_i w_{i+1}}
    %\]
    %where $\overline{w_i w_{i+1}}$ denotes the straight line segment connecting $w_i$ and $w_{i+1}$
    \[
    \gamma_{\rho_G}(a^*,b^*) = \{w_i\}^n_{i=1}.
    \]
    \end{nonsec}
%\end{defn}

%\textbf{Exact Curve (on G)}: Firstly we choose two points $a, b$ on $\mathbb{H}^2$, then with analytic formula, we obtain the hyperbolic geodesic $\gamma_{\rho_{\mathbb{H}^2}}(a,b) \subset \mathbb{H}^2$. After discretization(discretization needs to be clarified), we obtained a finite sequence $\{z_i\}_{i=1}^n \subset \gamma_{\rho_{\mathbb{H}^2}}(a,b)$, where $z_1 = a,\; z_n=b$. Then apply Schwarz-Christoffel Transformation $f: \mathbb{H}^2 \to G$ on $\{z_i\}_{i=1}^n $, to obtain $\{f(z_i)\}_{i=1}^n$ on $G$, where $f(z_1)=f(a)=a^*$ and $f(z_n)=f(b)=b^*$. Finally we obtain the hyperbolic geodesic $\gamma_{\rho_G}(a^*,b^*)$ on $G$ by joining the points $z_i$ in finite complex sequence $\{z_i\}_{i=1}^n$ with segments.

\begin{nonsec}{\bf Approximated hyperbolic geodesics.} \label{approximated_hyperbolic_geodesics}
 For any $a^*, b^* \in G$, the approximated hyperbolic geodesic $\gamma^*_{\rho_G}(a^*,b^*)$ on $G$ is constructed via the following discrete approximation scheme:
Let $a = f^{-1}(a^*)$, $b = f^{-1}(b^*)$. We define the distance between $z_i$ and $z_j$ 
%in the sense of Algorithm \ref{alg:combined_algorithm} 
as
\(
\rho_{\mathbb{H}^2}(f^{-1}(z_i),f^{-1}(z_j))
\)
 This distance is the weight function in the sense of Algorithm \ref{alg:combined_algorithm}. We apply Algorithm \ref{alg:combined_algorithm} on the discretized domain $G$ to obtain the discrete approximation of the length-minimizing curve $\gamma^*_{\rho_G}(a^*,b^*)$.
\end{nonsec}

%\begin{defn}[Approximated Hyperbolic Geodesic]
% For any $a^*, b^* \in G$, the approximated hyperbolic geodesic $\gamma^*_{\rho_G}(a^*,b^*)$ on $G$ is constructed via the following discrete approximation scheme:

%Let $a = f^{-1}(a^*)$, $b = f^{-1}(b^*)$. Then:
%\begin{itemize}
%    \item \textbf{Weight Function Definition}: Define the edge weight function for the graph as:
%    \[
%    w(z_i,z_j) = \rho_{\mathbb{H}^2}(f^{-1}(z_i),f^{-1}(z_j))
%    \]
%    where $\rho_{\mathbb{H}^2}$ is the hyperbolic distance metric on $\mathbb{H}^2$
    
 %   \item \textbf{Path Computation}: Apply \textbf{Algo 1, Algo 2} on the discretized domain $G$ to obtain $\gamma^*_{\rho_G}(a^*,b^*) = shortestpath$.
    
    %\item \textbf{Curve Construction}: The computed path provides the discrete approximation of the hyperbolic geodesic $\gamma_{\rho_G}(a^*,b^*)$ on $G$
%\end{itemize}
%\end{defn}

\begin{nonsec}{\bf Approximated quasihyperbolic geodesics.}
For any $a^*, b^* \in G$, the approximated quasihyperbolic geodesic $\gamma^*_{k_G}(a^*,b^*)$ on $G$ is constructed via Algorithm \ref{alg:combined_algorithm} where the distance between $z_i$ and $z_j$ is as in \eqref{qhwt}, and it is given by:
\[
\frac{|z_i - z_j|}{\min\{d(z_i,\partial G), d(z_j,\partial G)\}} \,.
\]
%in the sense of Algorithm \ref{alg:combined_algorithm}
This distance is the weight function in the sense of Algorithm \ref{alg:combined_algorithm}.
\end{nonsec}
%\textbf{QH Curve (on G)}: This means the approximated quasihyperbolic curve on G, similar to \textbf{DA Curve}, with the following weight function instead (induced from Def. of quasihyperbolic metric):
%\[
%w(z_i,z_j) = \frac{|z_i-z_j|}{\min\{d(z_i,\partial G), d(z_j,\partial G)\}}
%\]

%%%%%%%%%%%%%%%%%%%%%%
%%%%%%%%%%%%%%%%%%%%%%
%%%%%%%%%%%%%%%%%%%%%%
\section{Experiments on the Poincar\'e Unit Disk and Schwarz-Christoffel Mappings}
%%%%%%%%%%%%%%%%%%%%%%
%%%%%%%%%%%%%%%%%%%%%%
%%%%%%%%%%%%%%%%%%%%%%

\begin{nonsec}{\bf The unit disk.}
We first apply our algorithm to
approximate the quasihyperbolic metric of the unit disk $\BB$
equipped with the hyperbolic metric. In this case, it is known that for all $x,y \in \BB$   
\cite[Fig. 4.5., Remark 5.4]{hkv}: 
\[
\rho_{\BB}(x,y) \leq 2k_{\BB}(x,y) \leq 2\rho_{\BB} (x,y) \,.
\]
%For all $x,y$ on same diameter,
%\[
%\gamma_{\rho_{\BB}}(x,y)=\gamma_{k_{\BB}}(x,y)
%\]
%Otherwise,

The hyperbolic geodesics will be the sub-arcs of circular arcs orthogonal  to the boundary of $\BB$, joining the points $x$ and $y\,.$

%The latter is from the book "conformally invariant metrics and quasiconformal mapping" Remark 5.4.

%With designing several experiments, we obtain the following results: \textbf{Pictures (2 kinds geodesics)}.

Our aim is to design an experiment to show the convergence of Algorithm \ref{alg:combined_algorithm} on $\BB$ with decreasing $h$ and increasing $m$. With setting $p_1 = 0.9+0i$ and $p_2 = 0.495+0.495 i$, we obtain $\rho_{\BB}(p_1,p_2) \approx 2.9357$ by \eqref{rhoB}. Then we define $\epsilon_{h,m}(p_1,p_2) = |\rho_{\BB}(p_1,p_2) - \rho^*_{\BB}(p_1,p_2)|$ to obtain the error estimate between analytic results and approximated results.
\end{nonsec}
%\begin{table}[H]
%\begin{tabular}{@{}cccccc@{}}
%\toprule
%$h$      & $m$  & $\rho^*_{\BB}(p_1,p_2)$          & $k^*_{\BB}(p_1,p_2)$ & $\epsilon_{h,m}(p_1,p_2) $               & $\rho^*_{\BB}(p_1,p_2) / k^*_{\BB}(p_1,p_2)$        \\ \midrule
%0.01   & 3  & 2.94247           & 2.52575         & -0.00677              & 1.162308225 \\
%0.01   & 6  & 2.93933           & 2.52301         & -0.00363              & 1.163570497 \\
%0.01   & 8  & 2.94046           & 2.52355         & -0.00476              & 1.163321511 \\
%0.005  & 3  & 2.94162           & 2.52464         & -0.00592              & 1.162819253 \\
%0.005  & 6  & 2.93673           & 2.52100         & 0.001069391620290     & 1.164498215 \\
%0.005  & 8  & 2.93699           & 2.52111         & 0.001322206795484     & 1.164447406 \\
%0.005  & 10 & 2.93699           & 2.52111         & 0.001698595409156     & 1.164447406 \\
%0.0025 & 6  & 2.93654           & 2.52067         & 8.736021842490338e-04 & 1.164650668 \\
%0.0025 & 8  & 2.936170720955789 & 2.52042         & 5.064920359010294e-04 & 1.16476619  \\ \bottomrule
%\end{tabular}
%\caption{Result}
%\end{table}

\begin{table}[H]
\centering
\begin{tabular}{@{}cccccc@{}}
\toprule
$h$      & $m$  & $\rho^*_{\BB}(p_1,p_2)$ & $k^*_{\BB}(p_1,p_2)$ & $\epsilon_{h,m}(p_1,p_2)$ & $\rho^*_{\BB}(p_1,p_2) / k^*_{\BB}(p_1,p_2)$ \\ \midrule
0.01   & 3  & 2.9425  & 2.5258  & $6.8054\times 10^{-3}$ & 1.1623  \\
0.01   & 6  & 2.9393  & 2.5230  & $3.6631\times 10^{-3}$ & 1.1636  \\
0.01   & 8  & 2.9405  & 2.5236  & $4.7953\times 10^{-3}$ & 1.1633  \\
0.005  & 3  & 2.9416  & 2.5246  & $5.9590\times 10^{-3}$ & 1.1628  \\
0.005  & 6  & 2.9367  & 2.5210  & $1.0694\times 10^{-3}$  & 1.1645  \\
0.005  & 8  & 2.9370  & 2.5211  & $1.3222\times 10^{-3}$  & 1.1644  \\
0.005  & 10 & 2.9370  & 2.5211  & $1.6986\times 10^{-3}$  & 1.1644  \\
0.0025 & 6  & 2.9365  & 2.5207  & $8.7360\times 10^{-4}$  & 1.1647  \\
0.0025 & 8  & 2.9362  & 2.5204  & $5.0649\times 10^{-4}$  & 1.1648  \\ \bottomrule
\end{tabular}
\caption{Experimental results of local approximation of hyperbolic and quasihyperbolic geodesics and lengths in the unit disk. Approximated values and difference to the exact value for the hyperbolic case, are given.}
\label{hyp-disk}
\end{table}

In Table \ref{hyp-disk}, the errors in the measured length difference between the approximate and exact values decrease as $h$ decreases and $m$ increases. The accuracy reaches a reasonable level of $10^{-4}$. Another conclusion is that if we fix $h$ and only increase $m$, the improvement in accuracy may be reverted in some cases.

%%%%%%%%%%%%%%%%%%%%%%%%%%%%
%%%%%%%%%%%%%%%%%%%%%%%%%%%%
%%%%%%%%%%%%%%%%%%%%%%%%%%%%
\begin{nonsec}{\bf Polygonal domains.}
%%%%%%%%%%%%%%%%%%%%%%%%%%%%
%%%%%%%%%%%%%%%%%%%%%%%%%%%%
%%%%%%%%%%%%%%%%%%%%%%%%%%%%
Next, we consider two conformal mappings of the upper half-plane onto
a polygonal domain. By the Schwarz-Christoffel formula, such mappings exist and efficient methods for finding these conformal mappings are given in \cite{dritre}.

In our next two examples, we will apply the idea described in Subsection \ref{approximated_hyperbolic_geodesics} to the specific cases of polygonal domains. In these cases, the required conformal mapping can be obtained from the Schwarz-Christoffel formula. We use the toolbox of Driscoll \cite{Dri} to compute the mapping. This toolbox has a very high accuracy, so its effect on the computational errors is negligible.

Our first example is a polygonal quadrilateral and in this case  the formula is given in Theorem \ref{ellsintheo}.
Our second example is motivated by  Kythe \cite[pp. 226--227]{ky}. For this second example,
let $k  \in (0,1) ,$  $\xi_1 =1/k,  \xi_2  \in (0,1)$ and let
 \[
 A=   -\xi_1  , B=\-1/k, C=1,\text{ and } D=   -\xi_2,
  \]
and consider the function
\begin{equation}  \label{Kfun}
  f( \zeta)= -i  \zeta  \int_0^1  \frac{g ( \zeta,t,  \xi_2)  \qquad dt }{ g ( \zeta,t,  \xi_2) g ( \zeta,t,  \xi_1)  \sqrt{-(k^2) \zeta^2 t^2 +1}  }  \,,  \quad  g ( \zeta,t,  u) =
  \sqrt{ \zeta^2 t^2 + u^2  }.
\end{equation}
This function maps the upper half-plane onto a polygon, such that the points $A$, $B$, $C$, $D$, $-D$, $-C$, $-B$, $-A$ on the real axis are mapped onto the vertices $A'$, $B'$, $C'$, $D'$, $-D'$, $-C'$, $-B'$, $-A'$ of the polygon, see Figures \ref{figKy2} and \ref{figKy1}.
\end{nonsec}

\begin{figure}[H]
\includegraphics[scale=0.1,height=0.4\textheight]{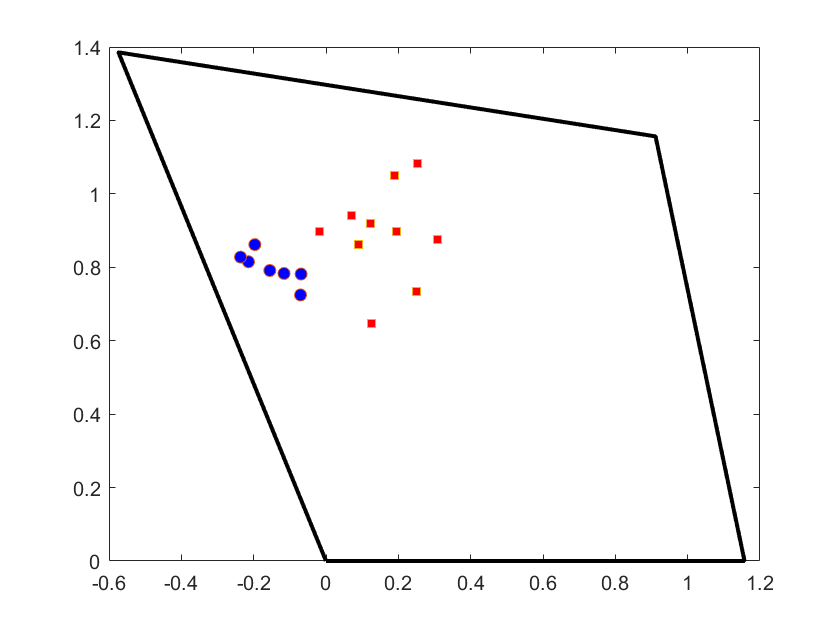} 
\caption{The quadrilateral of Example \ref{ex20250826}. }
\label{figHVV}
\end{figure}

\begin{nonsec}{\bf Example \cite{hvv}.} \label{ex20250826}
According to Theorem \ref{ellsintheo}, one can
map the upper half-plane onto a convex polygonal
quadrilateral by means of the formula \eqref{eq:Fint}. This transformation maps the
points $0,1, 1/r^2, \infty$ such that the points
$0$ and $1$ are mapped on the points $0$ and ${}_2F_1(a,b;c;r^2)$ of the real axis while the points $ 1/r^2, \infty$ are mapped into the upper half-plane. Observe the parameter 
constraints for the mapping  \eqref{eq:Fint}:
\[
0<a,b <1, 
\max\{a+b,1\}<c \le 1 + \min\{a,b\} \,.
\]
The quadrilateral illustrated in Figure \ref{figHVV} is convex, with the interior
angles equal to 
\[
b \pi, (c-b) \pi, (1-a) \pi,\text{ and } (1+a-c) \pi \,.
\]
\end{nonsec}

%\begin{table}[]
%\begin{tabular}{@{}cccccc@{}}
%\toprule
%$z_1$    & $z_2$   & $\rho_{\mathbb{H}^2}(z_1,z_2)$ & $f(z_1)$           & $f(z_2)$         & DM{[}f{[}z1{]},f{[}z2{]}{]} \\  \midrule
%-2+4i & 2i   & 0.962424        & -0.0689481+0.781559i  & 0.126475+0.648131i   & 0.96336                     \\
%-7+4i & 4+2i & 2.86608         & -0.236362+0.82775i    & 0.252413+1.08393i    & 2.86969                     \\
%-6+6i & 2+2i & 2.1459          & -0.196757+0.86178i    & 0.308001+0.874842i   & 2.14825                     \\
%-6+4i & 2+4i & 1.76275         & -0.214293+0.815009i   & 0.12258+0.92061i     & 1.76421                     \\
%-4+4i & 2+3i & 1.58481         & -0.155392+0.791223i   & 0.196+0.896825i      & 1.58567                     \\
%-2+4i & 1+2i & 1.20122         & -0.0689481+0.781559i  & 0.249751+0.735099i   & 1.20383                     \\
%-2+3i & 4+3i & 1.76275         & 5-0.0704523+0.724596i & 0.188413+1.05103i    & 1.76481                     \\
%-4+4i & 2+5i & 1.27253         & -0.155392+0.791223i   & 0.0700409+0.941958i  & 1.27541                     \\
%-7+4i & 6i   & 1.37511         & -0.236362+0.82775i    & -0.0188412+0.898678i & 1.37807                     \\
%-3+4i & 1+4i & 0.962424        & -0.116076+0.78314i    & 0.0899137+0.86272i   & 0.96369                     \\ \bottomrule
%\end{tabular}
%\caption{Example20250827, with $h=0.005,m=4$}
%\end{table}

After applying Algorithm \ref{alg:combined_algorithm} locally with setting $h=0.005,m=4$, we obtain the results below with defining $\epsilon=|\rho_{\UH}(z_1,z_2)-\rho_G^*\big(f(z_1),f(z_2)\big)|$ to describe the difference between the hyperbolic metric and the approximated hyperbolic metric.

\begin{table}[H]
\centering
\begin{tabular}{@{}ccccc@{}}
\toprule
$z_1$    & $z_2$   & $\rho_{\UH}(z_1,z_2)$ & $\rho_G^*\big(f(z_1),f(z_2)\big)$ & $\epsilon$ \\ \midrule
$-2+4i$ & $2i$    & 0.96242 & 0.96336 & $9.3957\times 10^{-4}$ \\
$-7+4i$ & $4+2i$  & 2.8661  & 2.8697  & $3.6090\times 10^{-3}$ \\
$-6+6i$ & $2+2i$  & 2.1459  & 2.1483  & $2.3568\times 10^{-3}$ \\
$-6+4i$ & $2+4i$  & 1.7628  & 1.7642  & $1.4674\times 10^{-3}$ \\
$-4+4i$ & $2+3i$  & 1.5848  & 1.5857  & $8.5910\times 10^{-4}$ \\
$-2+4i$ & $1+2i$  & 1.2012  & 1.2038  & $2.6128\times 10^{-3}$ \\
$-2+3i$ & $4+3i$  & 1.7628  & 1.7648  & $2.0600\times 10^{-3}$ \\
$-4+4i$ & $2+5i$  & 1.2725  & 1.2754  & $2.8754\times 10^{-3}$ \\
$-7+4i$ & $6i$    & 1.3751  & 1.3781  & $2.9542\times 10^{-3}$ \\
$-3+4i$ & $1+4i$  & 0.96242 & 0.96369 & $1.2620\times 10^{-3}$ \\ \bottomrule
\end{tabular}
\caption{Example \ref{ex20250826}.}
\end{table}

\begin{nonsec}{\bf Example, Kythe \cite[pp. 226--227]{ky}.} \label{ex20250325}
We construct a conformal map of the upper half-plane onto a given domain of Figure \ref{figKy2} so that the corresponding boundary points are as in Figure \ref{figKy2}. We consider 10 pairs of points $(z_j, w_j)$, $j=1,2,\ldots,10$, in the upper half-plane with integer coordinates and compute their hyperbolic distance. Then we use the
Kythe map $h$ from the upper half-plane (see Figure \ref{figKy1}) onto the domain of Figure \ref{figKy2}. After applying Algorithm \ref{alg:combined_algorithm} locally with the same arguments as in the previous example, we obtain the following results, with parameter values $\xi_1=2.0, \xi_2=0.7, k=0.8.$ The results are tabulated in Table \ref{table:ex_43}.

\begin{figure}
    \centering
		\begin{tikzpicture}[scale=2.5]
			%--BASIC SHAPE---------------------
			%------DEFINITIONS-----------------
			\coordinate (o) at (0,0);
			\coordinate (d) at (0,0.2212);
			\coordinate (c) at (-0.3703,0.2212);
			\coordinate (B) at (-0.3630,1.1490);
			\coordinate (A) at (1.0892,1.1448);
			\coordinate (md) at (0,-0.2212);
			\coordinate (mc) at (-0.3703,-0.2212);
			\coordinate (mB) at (-0.3630,-1.1490);
			\coordinate (mA) at (1.0892,-1.1448);
			%-------LINES-----------------------
			\draw[step=0.5,gray,thin] ($(B)+(-.2,.2)$) grid ($(mA)+(.2,-.2)$);
			\foreach \y in {-1,-0.5,0,0.5,1} \draw (-.65,\y) node[] {\tiny $\y$};
			\foreach \x in {-0.5,0,0.5,1} \draw (\x,-1.4) node[] {\tiny $\x$};
			\draw[thick] (o) -- (d) -- (c) -- (B) -- (A) -- (mA) -- (mB) -- (mc) -- (md)-- (o);
			%-------NODES-----------------------
			\draw \foreach \x in {(o),(d),(c),(B),(A),(md),(mc),(mB),(mA)} {
				\x node[circle, fill=black, inner sep=0pt, minimum width=5pt] {} }; % black nodes 
			%--END OF BASIC SHAPE---------------
			%-----------------------------------
			%--LABELS--------------------
			\node at (o) [left = 1mm of o] {$O$};
			\node at (d) [above = 1mm of d] {$D'$};
			\node at (c) [left = 1mm of c] {$C'$};
			\node at (B) [left = 1mm of B] {$B'$};
			\node at (A) [right = 1mm of A] {$A'$};
			\node at (md) [below = 1mm of md] {$-D'$};
			\node at (mc) [left = 1mm of mc] {$-C'$};
			\node at (mB) [left = 1mm of mB] {$-B'$};
			\node at (mA) [right = 1mm of mA] {$-A'$};
			%--LENGTH LABELS AND BRACES---------
			%--ANGLES---------------------------
		\end{tikzpicture}
\caption{This graphics shows the figure from Kythe \cite[pp. 226--227]{ky}. The labels $A',B',\dots$ indicate the image points under the mapping \eqref{Kfun} of the pre-vertex points  $A,B,\dots$ on the real axis. }
\label{figKy2}
\end{figure}
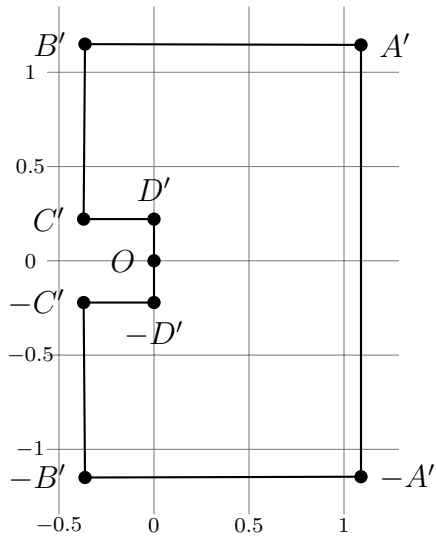 

%%%%%%%%%%%%%%%%%%%%%%%%%%%%

%\textbf{Replace the following as Tikz}
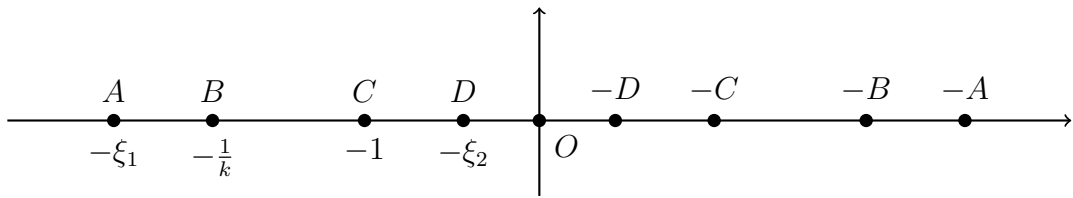
\begin{figure}[H]
\begin{tikzpicture}
    \coordinate (o) at (0,0);
    \coordinate (d) at (-1,0);
    \coordinate (md) at ($-1*(d)$);
    \coordinate (c) at ($(d)+(-1.3,0)$);
    \coordinate (mc) at ($-1*(c)$);
    \coordinate (B) at ($(c)+(-2,0)$);
    \coordinate (mB) at ($-1*(B)$);
    \coordinate (A) at ($(B)+(-1.3,0)$);
    \coordinate (mA) at ($-1*(A)$);
    \draw[->, thick] (-7,0) -- (7,0);
    \draw[->, thick] (0,-1) -- (0,1.5);
    %-------NODES-----------------------
	\draw \foreach \x in {(o),(d),(md),(c),(mc),(B),(mB),(A),(mA)} {
			\x node[circle, fill=black, inner sep=0pt, minimum width=5pt] {} }; % black nodes 
    %-------LABELS----------------------
	\node at (o) [below right = .7mm of o] {$O$};
	\node at (d) [above = 1mm of d] {$D$};
	\node at (md) [above = 1mm of md] {$-D$};
	\node at (c) [above = 1mm of c] {$C$};
	\node at (mc) [above = 1mm of mc] {$-C$};
	\node at (B) [above = 1mm of B] {$B$};
	\node at (mB) [above = 1mm of mB] {$-B$};
	\node at (A) [above = 1mm of A] {$A$};
	\node at (mA) [above = 1mm of mA] {$-A$};
	\node at (d) [below = 1mm of d] {$-\xi_2$};
	\node at (c) [below = 1mm of c] {$-1$};
	\node at (B) [below = 1mm of B] {$-\frac{1}{k}$};
	\node at (A) [below = 1mm of A] {$-\xi_1$};
\end{tikzpicture}
\caption{This graphics reproduces the figure from Kythe \cite[p. 227]{ky}. The pre-vertices  $A,B,\dots$ are on the real axis. }
\label{figKy1}
\end{figure}

%\begin{table}[htbp]
%\centering
%\small
%\begin{tabular}{cccccc}
%\toprule
%$z_1$    & $z_2$   & $\rho_{\mathbb{H}^2}(z_1,z_2)$ & $f(z_1)$           & $f(z_2)$         & DM{[}f{[}z1{]},f{[}z2{]}{]} \\ 
%\midrule
%-1+3i & 4+3i & 1.5170          & 0.73708 + 0.098334i & 0.93043 - 0.19568i & 1.5316                      \\
%-1+4i & 3+3i & 1.1293          & 0.80836 + 0.063047i & 0.87107 - 0.19558i & 1.1427                      \\
%-1+4i & 2+2i & 1.2012          & 0.80836 + 0.063047i & 0.74855 - 0.26696i & 1.2132                      \\
%-2+4i & 4+2i & 1.9248          & 0.84306 + 0.11205i  & 0.94993 - 0.25148i & 1.9545                      \\
%4i    & 2+2i & 0.96242         & 0.79503             & 0.74855 - 0.26696i & 0.9768                      \\
% -2+4i & 4+2i & 1.9248          & 0.84306 + 0.11205i  & 0.94993 - 0.25148i & 1.9545                      \\
%-2+4i & 2+2i & 1.4506          & 0.84306 + 0.11205i  & 0.74855 - 0.26696i & 1.4703                      \\
%2i    & 2+2i & 0.96242         & 0.57552             & 0.74855 - 0.26696i & 0.9693                      \\
%-1+2i & 3+2i & 1.7627          & 0.62469 + 0.16678i  & 0.87107 - 0.19558i & 1.7764                      \\
%-2+4i & 4+2i & 1.9248          & 0.84306 + 0.11205i  & 0.94993 - 0.25148i & 1.9545                      \\ \bottomrule
%\end{tabular}
%\caption{Example20250325}
%\end{table}

\begin{table}[H]
\begin{tabular}{@{}ccccc@{}}
\toprule
$z_1$    & $z_2$   & $\rho_{\UH}(z_1,z_2)$ & $\rho_G^*\big(f(z_1),f(z_2)\big)$ & $\epsilon$ \\ \midrule
$-1+3i$ & $4+3i$ & 1.5170 & 1.5316 & $1.0155\times 10^{-3}$ \\
$-1+4i$ & $3+3i$ & 1.1293 & 1.1307 & $1.4684\times 10^{-3}$ \\
$-1+4i$ & $2+2i$ & 1.2012 & 1.2025 & $1.2997\times 10^{-3}$ \\
$4i$    & $2+2i$ & 0.96242 & 0.96355 & $1.1256\times 10^{-3}$ \\
$-2+4i$ & $2+2i$ & 1.4506 & 1.4522 & $2.1042\times 10^{-3}$ \\
$2i$    & $2+2i$ & 0.96242 & 0.96353 & $1.1020\times 10^{-3}$ \\
$-1+2i$ & $3+2i$ & 1.7627 & 1.7645 & $1.7569\times 10^{-3}$ \\
$-2+4i$ & $4+2i$ & 1.9248 & 1.9545 & $1.7300\times 10^{-3}$  \\ \bottomrule
\end{tabular}
\caption{Example \ref{ex20250325}, with $h=0.005,m=4.$}\label{table:ex_43}
\end{table}

%%%%%%%%%%%%%%%%%%%%%%%%%%%%
%%%%%%%%%%%%%%%%%%%%%%%%%%%%

%\textbf{Here if you want, we can design a table with changing $h,m$ and observing how $\epsilon$ changes.}
%%%%%%%%%%%%%%%%%%%%%%%%%%%%
%%%%%%%%%%%%%%%%%%%%%%%%%%%%
%%%%%%%%%%%%%%%%%%%%%%%%%%%%

\section{Bifurcation and Non-uniqueness of Geodesics}
%%%%%%%%%%%%%%%%%%%%%%%%%%%%
%%%%%%%%%%%%%%%%%%%%%%%%%%%%
%%%%%%%%%%%%%%%%%%%%%%%%%%%%
%\textbf{Def. (Theoretical Bifurcation Point) Needed to be completed.}

An interesting phenomenon arising from quasihyperbolic geometry is bifurcation of geodesics. This may arise in two specific situations:
\begin{enumerate}
\item[(a)] If the domain has more than one boundary component, it is possible that there are more than one quasihyperbolic geodesic between individual points choosing different routes around interior boundary components. In this case, the geodesics may go to different directions from the start, or they may have a common part and bifurcate later.
\item[(b)] It is possible that even on a simply connected domain $G$ quasihyperbolic geodesics connecting a point $x$ to $y$ and one connecting it to a different point $z$ coincide at the points of some sub-curve, and only then separate their ways.
\end{enumerate}

Non-uniqueness is possible for hyperbolic geodesics in the case if the domain is multiply connected, i.e. if the domain has  more than one boundary component. Note that hyperbolic geodesics joining $x$ to separate points $y,z$ may never have a common sub-curve. For hyperbolic domains conformally equivalent to the disk, this is easy to see, because conformal mappings are homeomorphisms, and the bifurcation property does not arise for the case of the disk, where explicit form of the geodesics is known. But the situation is very different in the case of quasihyperbolic geodesics, where bifurcation is possible already in the case of convex domains.

In this section, our goal is to experimentally demonstrate bifurcation of geodesics in several cases that are known from the literature. The examples are mostly reproductions of theoretical results from \cite{AnJa} and \cite{Linden}. We will start with a formal definition.

\begin{nonsec}{\bf Definition.}\label{bifurcation point}
Let $\gamma_1, \gamma_2$ be regular curves on a domain $G$, so that $\gamma_1,\gamma_2:[0,1] \to G$. Suppose that there exists a point $p \in G$ s.t. $\gamma_1(t_1)=\gamma_2(t_2)=p$ for some $t_1,t_2 \in (0,1)$, the point $p$ is a  {\it bifurcation point} for $\gamma_1$,$\gamma_2$ if there exists $\epsilon >0$ and a smooth strictly increasing function $\varphi: [t_1-\epsilon,t_1+\epsilon]\to [0,1]$ with $\varphi(t_1)=t_2$ so that either 
\[
\left\{\begin{array}{rcl}
\gamma_1'(t) = \gamma_2'\big(\varphi(t)\big), &\textrm{for}& t\in(t_1-\epsilon,t_1), \\
\gamma_1'(t) \neq \gamma_2'\big(\varphi(t)\big), &\textrm{for}& t\in(t_1,t_1+\epsilon),
\end{array}\right.
\]
or
\[
\left\{\begin{array}{rcl}
\gamma_1'(t) \neq \gamma_2'\big(\varphi(t)\big), &\textrm{for}& t\in(t_1-\epsilon,t_1), \\
\gamma_1'(t) = \gamma_2'\big(\varphi(t)\big), &\textrm{for}& t\in(t_1,t_1+\epsilon).
\end{array}\right.
\]
If $t_1$ or $t_2$ is equal to $0$ or $1$, we may extend this definition so that only the inequality condition is required.
\end{nonsec}

Based on the definition above, we design the following method to approximate the bifurcation point. It is noted that geodesics of quasihyperbolic metric are smooth, so the above definition can be applied. In the multiply connected case it is possible that geodesics between points $x$ and $y$ are not unique, and they may touch at individual points that are not considered as bifurcation points in the sense of this definition.

It should be noted that hyperbolic geodesics may not bifurcate on simply connected domains. This follows from the properties of hyperbolic geodesics on the unit disk combined with their conformal invariance.

%\textbf{Def. (approximated Bifurcation Point)} 
%Given $\gamma_1 = \{z_i\}_{i=1}^m$ and $\gamma_2 = \{w_j\}_{j=1}^n$ generated from DM, with corresponding argument sequences $\{k_{i,\gamma_1} \in [-\pi,\pi]: k_{i,\gamma_1}  = \arg(z_{i+1} - z_i),z_i \in \gamma_1 \}_{i=1}^{m-1}$ and $\{k_{j,\gamma_2} \in [-\pi,\pi]: k_{j,\gamma_2}  = \arg(w_{j+1} - w_j),w_j \in \gamma_2 \}_{j=1}^{n-1}$. Given a point $p \in \gamma_1 \cup\gamma_2 $ s.t. $\exists i,j$ with $1\leq i \leq m-1,1 \leq j \leq n-1$ s.t. $p \in [z_i,z_{i+1}]$ and $p \in[w_j,w_{j+1}]$ (here [a,b] be segments). More precisely,$p \in [\operatorname{argmin}_{x \in \{z_i,w_j\}}|x-p|,\operatorname{argmin}_{x \in \{z_{i+1},w_{j+1}\}}|x-p|]=[a,b]$. The \textbf{approximated bifurcation point} should firstly hold $k_{i,\gamma_1} = k_{j,\gamma_2}$, then:
%\begin{enumerate}
%    \item If $k_{i-1,\gamma_1} = k_{j-1,\gamma_2}, k_{i+1,\gamma_1} \neq k_{j+1,\gamma_2}$, \textbf{approximated bifurcation point} is $b$.
%    \item If $k_{i-1,\gamma_1} \neq k_{j-1,\gamma_2}, k_{i+1,\gamma_1} = k_{j+1,\gamma_2}$, \textbf{approximated bifurcation point} is $a$.
%\end{enumerate}
\small

\begin{algorithm}[H]
\caption{Approximated Bifurcation Point Identification}
\label{alg:bifurcation_point}
\begin{algorithmic}[1]
\REQUIRE 
\STATE $\gamma_1^* = \{z_i\}_{i=1}^m$: first path generated by Dijkstra's algorithm
\STATE $\gamma_2^* = \{w_j\}_{j=1}^n$: second path generated from Dijkstra's algorithm  
\STATE $\{k_{i,\gamma_1^*} \in [-\pi,\pi]: k_{i,\gamma_1^*} = \arg(z_{i+1} - z_i), z_i \in \gamma_1^* \}_{i=1}^{m-1}$: argument sequence for $\gamma_1^*$
\STATE $\{k_{j,\gamma_2^*} \in [-\pi,\pi]: k_{j,\gamma_2^*} = \arg(w_{j+1} - w_j), w_j \in \gamma_2^* \}_{j=1}^{n-1}$: argument sequence for $\gamma_2^*$
\STATE $p \in \gamma_1^* \cup \gamma_2^*$: point satisfying $p \in [z_i,z_{i+1}]$ and $p \in [w_j,w_{j+1}]$ for some $1\leq i \leq m-1, 1 \leq j \leq n-1$
\ENSURE Approximated bifurcation point

\STATE $a \gets \operatorname{argmin}_{x \in \{z_i,w_j\}}|x-p|$
\STATE $b \gets \operatorname{argmin}_{x \in \{z_{i+1},w_{j+1}\}}|x-p|$
\STATE $[a,b]$ defines the segment containing $p$

\IF{$k_{i,\gamma_1^*} = k_{j,\gamma_2^*}$}
    \IF{$k_{i-1,\gamma_1^*} = k_{j-1,\gamma_2^*}$ \AND $k_{i+1,\gamma_1^*} \neq k_{j+1,\gamma_2^*}$}
        \STATE \textbf{return} $b$ \COMMENT{Case 1: approximated bifurcation point is $b$}
    \ELSIF{$k_{i-1,\gamma_1^*} \neq k_{j-1,\gamma_2^*}$ \AND $k_{i+1,\gamma_1^*} = k_{j+1,\gamma_2^*}$}
        \STATE \textbf{return} $a$ \COMMENT{Case 2: approximated bifurcation point is $a$}
    \ENDIF
\ENDIF

\STATE \textbf{return} None \COMMENT{No bifurcation point found under given conditions}
\end{algorithmic}
\end{algorithm}

\normalsize

\begin{nonsec}{\bf Remark.}
    Assume there exists a bifurcation point $p$ between two smooth curves $\gamma_1$ and $\gamma_2$ defined as in Definition \ref{bifurcation point}. When applying Algorithm \ref{alg:combined_algorithm} to obtain approximated geodesics $\gamma_1^*$ and $\gamma_2^*$, and the approximated bifurcation point derived from Algorithm \ref{alg:bifurcation_point} will converge to $p$ when refinements are done in Algorithm \ref{alg:combined_algorithm}.
\end{nonsec}

%\begin{nonsec}{\bf Multiply Connected Domains.}
%%%%%%%%%%%%%%%%%%%%%%%%%%%%
%%%%%%%%%%%%%%%%%%%%%%%%%%%%
%%%%%%%%%%%%%%%%%%%%%%%%%%%%
\begin{nonsec}{\bf Domains with Punctures.}
\label{sec: Lemma 5.1}
An application of our method is to demonstrate numerically the results originating from the proof of \cite[Lemma 5.1]{AnJa}. Let $G_n=R_n \setminus (P_n \cup L_n)$, where $R_n = \{a+bi \; | \; a\in(-1,(n-2)\sqrt3+1),b \in (-1,1)\}$, $P_n=\{0,\sqrt3,\dots,(n-2)\sqrt3\}$ and $L_n=\{a+bi \; | \; a \in \{ \frac {\sqrt 3} 2,\sqrt{3} + \frac {\sqrt 3} 2, \dots,(n-3)\sqrt{3} + \frac {\sqrt 3} 2 \},b\in(-1,-\frac 1 2) \cup (\frac 1 2 , 1)\}$.
Let $x = -\frac 1 2, y=(n-2)\sqrt 3+ \frac 1 2$, our aim is to approximate $\gamma_k(x,y)$ by applying Algorithm \ref{alg:combined_algorithm} on this domain.

% \begin{figure}[H]
%  \centering
%  \includegraphics[width=\linewidth]{figs/bifurcation/multiply connected domain/AnJa_fig3.png}
%  \caption{\cite[Figure 3. Geodesic $\gamma_1$ in the proof of Lemma 5.1]{AnJa}}
%  \label{fig: Lemma51}
%\end{figure}

Setting the values $h=0.025, m=4$, we obtain an approximation of the geodesic illustrated in Figure \ref{fig: Lemma51_DM}.
\end{nonsec}

\begin{figure}[H]
  \centering
  \includegraphics[width=0.5\linewidth]{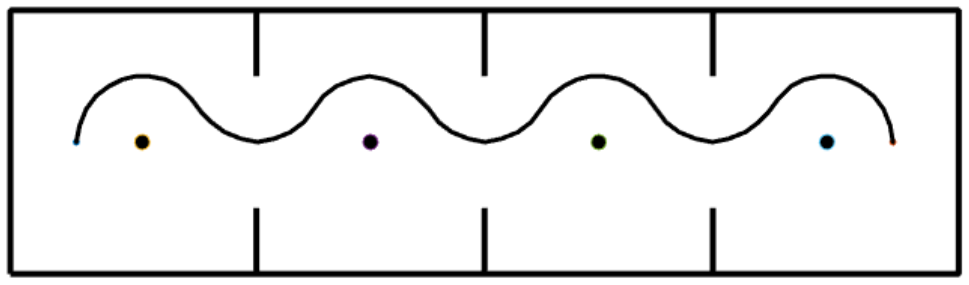}
  \caption{Algorithm \ref{alg:combined_algorithm} on \cite[Lemma 5.1]{AnJa}, with $h=0.025$ and $m=4$.}
  \label{fig: Lemma51_DM}
\end{figure}

Although the pictures look very similar to results in \cite[Lemma 5.1]{AnJa}, we still want to find a quantitave error estimate to compare our results to those in \cite{AnJa}. The following $\epsilon$ means the total error between $\{z_i\}_{i=1}^n$ and $\gamma_1$, where $\gamma_1$ defined in \cite[Lemma 5.1 Proof]{AnJa}. If $d(z_i, \gamma_1)$ is the distance between $z_i$ and $\gamma_1$, then we may use the error estimate from
\[
\epsilon = \sum_{i=1}^n d(z_i, \gamma_1).
\]
%\textbf{Here should be changed to data}
%\textbf{Add much more comments and details at here}
For the case $h=0.025, m=4$, the error is $\epsilon \approx  4.7577\cdot 10^{-5}$, which shows that the points of approximated quasihyperbolic geodesic are very close to $\gamma_1$.
%Results in this example will also support our conclusion in the following section. 
%\textbf{Add much more comments and details at here}

Based on \cite[Lemma 5.1]{AnJa}, with using the same domain $G$ as defined in \ref{sec: Lemma 5.1} and applying Algorithm \ref{alg:combined_algorithm} with $h=0.025$ and $m=4$, we obtain Figure \ref{fig:bif_multi_connect}.
\end{nonsec}

\begin{figure}[H]
  \centering
  \includegraphics[width=0.7\linewidth]{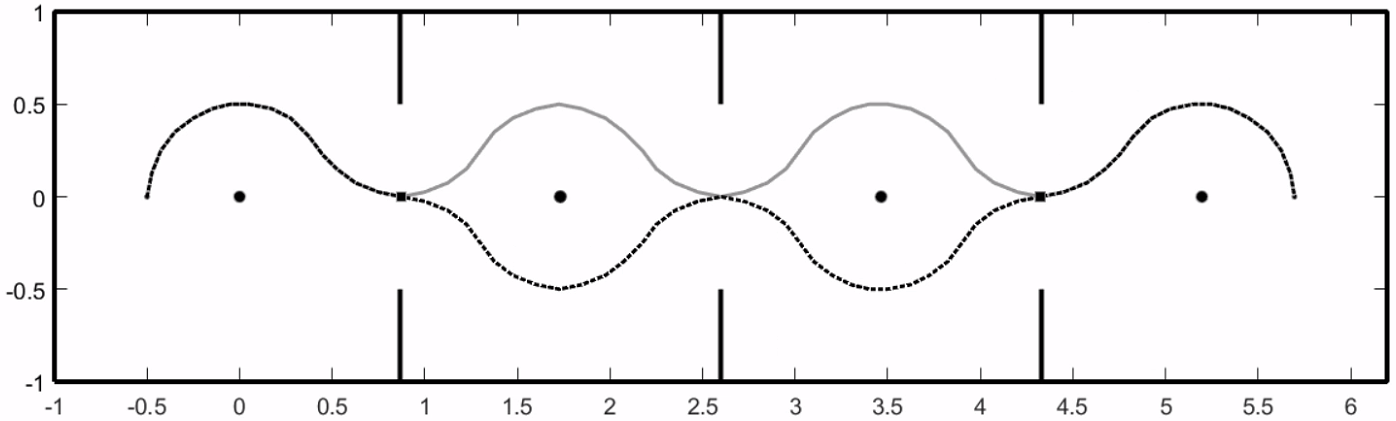}
  
  \caption{Non-uniqueness of geodesics on a multiply-connected domain. In this case, the geodesic may have two or more bifurcation points depending on the taken route.}\label{fig:bif_multi_connect}
\end{figure}

The theoretical bifurcation points appear at $\frac {\sqrt3} 2$ and $\frac {5\sqrt3} 2$ that are the same as in the proof of \cite[Lemma 5.1]{AnJa}. Moreover, the approximated bifurcation points are $0.875$ and $4.325$ as the square points shown in Figure \ref{fig:bif_multi_connect}, which respectively have errors $0.0089$ and $0.0051$.

\begin{figure}[H]
  \centering
  \begin{subfigure}{0.49\textwidth}
    \centering
    \includegraphics[width=\linewidth]{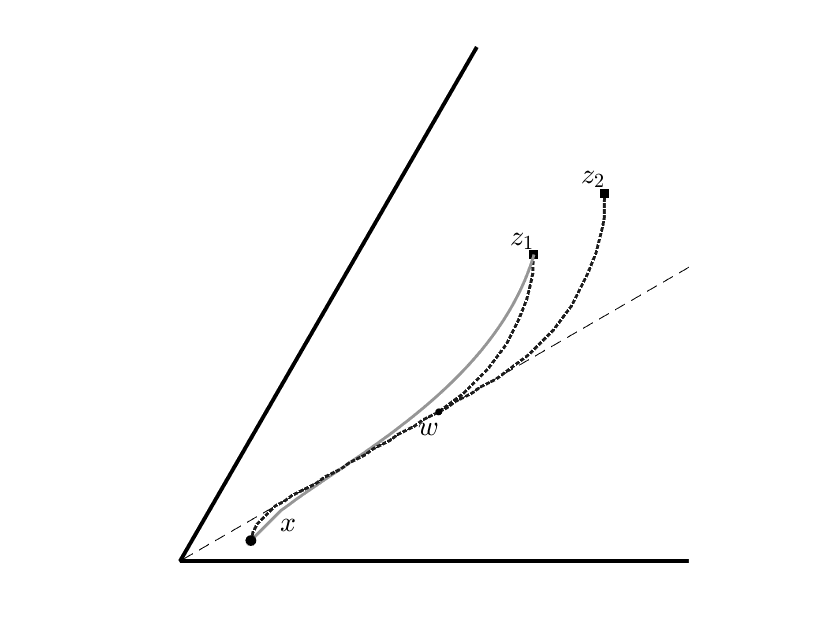}
    \caption{Bifurcation for points on the same side of bisector.}
    \label{fig: bif_conv_1}
  \end{subfigure}
  \hfill
  \begin{subfigure}{0.49\textwidth}
    \centering
    \includegraphics[width=\linewidth]{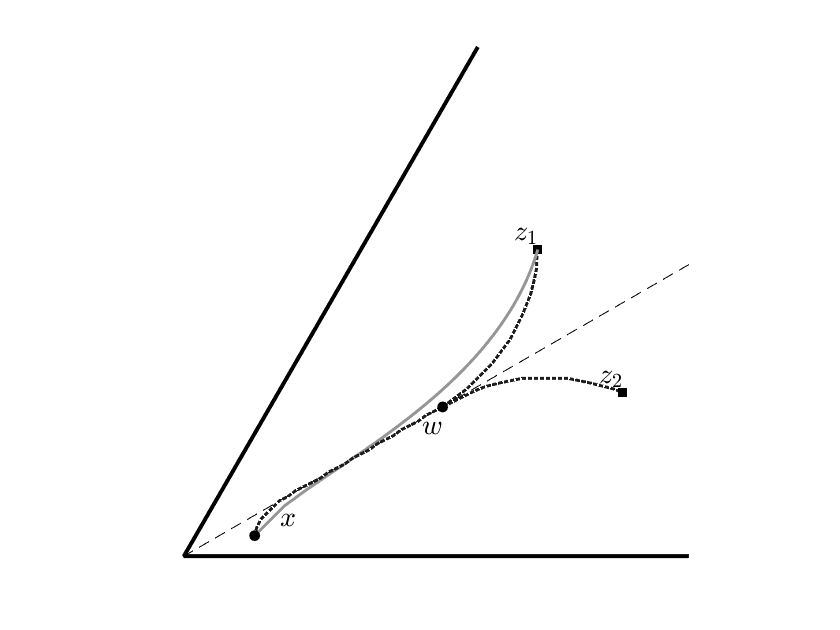}
    \caption{Bifurcation for points on different sides of bisector.}
    \label{fig: bif_conv_2}
  \end{subfigure}
  \caption{Bifurcation analysis in convex domain}
  \label{fig:bifurcation_convex}
\end{figure}

\begin{nonsec}{\bf Convex domains.}
    The idea of this subsection is from \cite{Linden}. The aim is to demonstrate that convexity of the domain does not guarantee uniqueness of geodesics. We consider an angular domain, with the angle between the sides $\frac \pi 3$. Let $x$ be the start point, let $w$ be the bifurcation point, and let $z_1$, $z_2$ be the end points, see Figure \ref{fig:bifurcation_convex}.
    
    The experiment applies Algorithm \ref{alg:combined_algorithm} locally with parameter values $h = 0.05$ and $m = 4$. The results are illustrated in Figure \ref{fig:bifurcation_convex}.  The solid line is the hyperbolic geodesic, and the dotted curve is the quasihyperbolic geodesic.

\end{nonsec}   

%\textbf{Need a Tikz here to replace the following image}

%\begin{figure}[H]
%    \centering
%    \includegraphics[width=0.5\linewidth]{figs/bifurcation/convex domain/angle_instruction.png}
%    \caption{Caption}
%    \label{fig:placeholder}
%\end{figure}
\begin{figure}[H]
    \centering
    \begin{tikzpicture}[scale=0.8]
		%--BASIC SHAPE---------------------
		%------DEFINITIONS-----------------
		\def\fii{50}; 
		\def\sidelength{11}; % length of the displayed domain sides
		\def\yl{6}; % distance of y from the origin
		\def\ra{{\yl*tan(\fii/2)}}; % radius of the circle
		\def\coeff{1/4}; % coefficient for theta
		\def\thet{\coeff*90+\coeff*\fii/2};
		\def\mthet{-\coeff*90-\coeff*\fii/2}; % minus theta, for some reason -\thet gives the wrong answer
		\coordinate (o) at (0:0); % origin
		\coordinate (s0) at (-\fii/2:\sidelength); % right side endpoint
		\coordinate (s1) at (\fii/2:\sidelength); % left side endpoint
		\coordinate (bis) at (0:\sidelength-1); % bisector endpoint
		\coordinate (x) at (0:2);
		\coordinate (y) at (0:\yl);
		\coordinate (a) at ($(y) + (-90:\ra)$); % centre of the circle on the right side of the area
		\coordinate (a') at ($(y) + (90:\ra)$); % centre of the circle on the left side of the area
		\coordinate (z0) at ($(a) + (-\fii/2:\ra)$);
		\coordinate (zj) at ($(a) + (-\fii/2+\thet:\ra)$);
		\coordinate (zjs) at ($(a') + (\fii/2+\mthet:\ra)$);
		\coordinate (zmid) at ($0.5*(zj) + 0.5*(zjs)$); % midpoint of the line segment between z_j and z*_j
		%-------LINES-----------------------
		\draw[thick] (s0) -- (o) -- (s1); % boundary
		\draw[thick, dashed] (o) -- (bis); % bisector
		\draw[thick, dotted] (z0) arc (-\fii/2:180-\fii/2:\ra); % semicircle
		\draw[thick, dashed] (zj) -- (zjs); % z_j to z_j*
		\draw[dashed] (y) -- (a) node[midway, right] {$r$}; % line from a to y
		\draw[dashed] (a) -- (zj)  node[midway, above] {$r$}; % line from a to z_j
		%-------NODES-----------------------
		\draw \foreach \x in {(o),(x),(y),(a),(z0),(zj),(zjs)} {
			\x node[circle, fill=black, inner sep=0pt, minimum width=5pt] {} }; % black nodes 
		%--END OF BASIC SHAPE---------------
		%-----------------------------------
		%--LABELS--------------------
		\draw (bis) node[right] {\emph{medial axis}};
		\node at (x) [below = 1mm of x] {$x$};
		\node at (y) [above right = .7mm of y] {$y$};
		\node at (z0) [below = 1mm of z0] {$z_0$};
		\node at (zj) [right = 1mm of zj] {$z_j$};
		\node at (zj) [right = 1mm of zjs] {$z_j^*$};
		\node at (o) [left = 1mm of o] {$O$};
		%--LENGTH LABELS AND BRACES---------
		\draw [thick, decorate, decoration = {calligraphic brace,raise=3pt,amplitude=6pt}] (o) -- (y) node[pos=0.5,above = 7pt,black] {$\ell$};
		\draw [thick, decorate, decoration = {calligraphic brace,raise=12pt,amplitude=6pt}] (a) -- (o) node[pos=0.4,below left=16pt,black] {$\frac{r}{\sin(\frac{\varphi}{2})}$};
		%--ANGLES---------------------------
		\draw pic [draw, angle radius=13 pt] {angle=s0--o--s1}; % angle phi
		\draw ($(o) + (0.2,-.4)$) node[] {$\varphi$}; %  phi label
		\draw pic [draw, angle radius=23 pt] {angle=s0--o--bis}; % angle phi/2
		\draw ($(o) + (0.6,-.6)$) node[] {$\frac{\varphi}{2}$}; % phi/2 label
		\draw pic [draw, angle radius = 6 pt] {right angle=o--y--a}; % right angle next to y
		\draw pic [draw, angle radius = 18 pt] {angle=z0--a--zj}; % theta_j
		\draw ($(a) + (1,-.21)$) node[] {$\theta_j$}; %  theta_j label
		\draw pic [draw, angle radius = 16 pt] {angle=y--a--o}; % y to center to origin
		\draw ($(a) + (-0.6,.7)$) node[] {$\frac{\pi-\varphi}{2}$};
		\draw pic [draw, angle radius = 6 pt] {right angle=y--zmid--zj}; % right angle with the medial axis and z_j--z*_j
	\end{tikzpicture}
    \caption{The angular domain.}
    \label{fig:angular}
\end{figure}
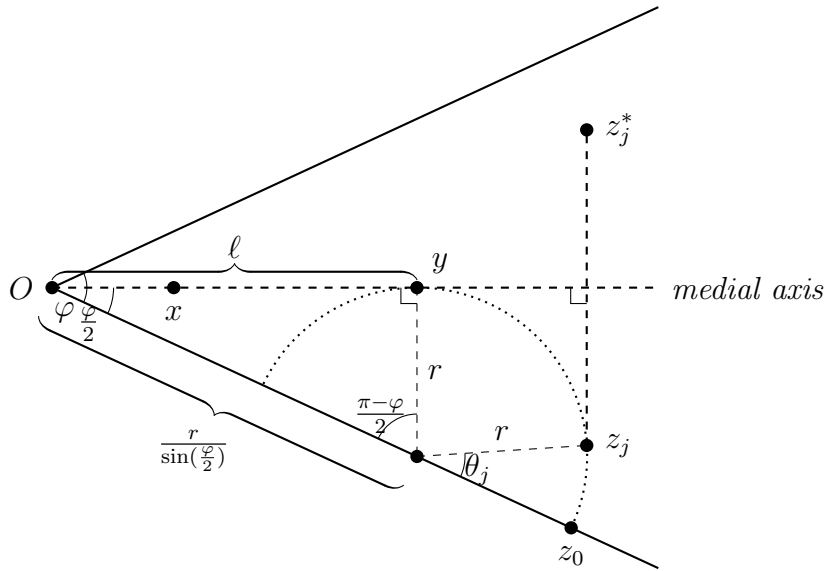

%\begin{figure}[H]
%  \centering
%  \includegraphics[width=0.5\linewidth]{figs/bifurcation/convex domain/angle_bifurcation_1_grayscale.png}
%  \caption{Bifurcation for points on the same side of bisector.}
%  \label{fig: bif_conv_1}
%\end{figure}

%\begin{figure}[H]
%  \centering
%  \includegraphics[width=0.5\linewidth]{figs/bifurcation/convex domain/angle_bifurcation_2_grayscale.png}
%  \caption{Bifurcation for points on different sides of bisector.}
%  \label{fig: bif_conv_2}
%\end{figure}

We next experimentally study the ideas from \cite[Theorem 4.6. Case 4]{Linden}. Consider for \( \varphi \in (0, \pi) \)  
and \( L > 0 \) a sequence of points in the sector  
$\{ z : |\arg z| < \frac{\varphi}{2},|z| < L \}$
constructed as follows  
\[
z_j = (\ell - ir) + re^{i(-\frac{\varphi}{2} + \theta_j)}, \quad 
\theta_j = \frac{j}{n} \cdot \frac{\pi + \varphi}{2}, \quad j = 1, \dots, n
\]  
where \( r = \ell \tan \frac{\varphi}{2},\ 0<\ell \ll L\) (See Figure \ref{fig:angular}).  

For the experiment, set \( L = 10^4,\ \ell = 3,\ n = 10 \),  
and \( \varphi = 2 \arctan\left(\frac{1}{2}\right) \). Also define \( z_j^* := \bar{z}_j \).  
The point \( y = \ell \) is the theoretical bifurcation  
point of the two geodesics joining  
\( x \) and \( z_j \), \( x \) and \( z_j^* \), respectively, when $x=\frac 1 2$.

Denote the approximated bifurcation points corresponding to the geodesic between $x$ and $z_j$ and $x$ and
$z^*_j$,  respectively, by  $w_j$. Our aim is to check the following:
\begin{enumerate}
    \item The change of $w_j$ when increasing $j$ from $1$ to $n$.
    \item Whether the trace of $\gamma_{x,z_j}$ and $\gamma_{x,z_j^*}$ agrees with the \cite[Theorem 4.6. Case 4]{Linden}.
\end{enumerate}

Define the error estimate $\epsilon_{\gamma_1^*}$ by the formula
\[
\epsilon_{\gamma_1^*} = \frac 1 {|S|}\sum_{p_j \in S} \bigg||p_j-(\ell+r \cdot i)|-r \bigg|,
\]
where $S=\{p_j \in \gamma_1^* \; | \; \Im(p_j) \neq 0\}$, and $\epsilon_{\gamma_2^*}$ is defined similarly. Both error estimates mean the average error between the theoretical quasihyperbolic geodesic and the approximated quasihyperbolic geodesic.

\begin{figure}[H]
  \centering
  \begin{subfigure}{0.48\textwidth}
    \centering
    \begin{tabular}{@{}cccc@{}}
    \toprule
    $j$ & $w_j$ & $\epsilon_{\gamma_1^*}$ & $\epsilon_{\gamma_2^*}$ \\ \midrule
    0 & 2.9        & $0.1049$        & $0.1049$        \\
    1 & 2.9        & $0.0342$        & $0.0342$        \\
    2 & 2.9        & $0.0276$        & $0.0276$        \\
    3 & 2.9        & $0.0304$        & $0.0304$        \\
    4 & 2.9        & $0.033$         & $0.033$         \\
    5 & 3.1        & $0.02$          & $0.02$          \\
    6 & 3.1        & $0.0245$        & $0.0245$        \\
    7 & 3.1        & $0.0071$        & $0.0071$        \\
    8 & 3.1        & $0.0042$        & $0.0042$        \\
    9 & 3.1        & $2.2204\times 10^{-16}$    & $2.2204\times 10^{-16}$   \\ \bottomrule
    \end{tabular}
    \caption{The case $h=0.025$ and $m=8$.}
    \label{tab:convex_params}
  \end{subfigure}
  \hfill
  \begin{subfigure}{0.48\textwidth}
    \centering
    \includegraphics[width=\linewidth]{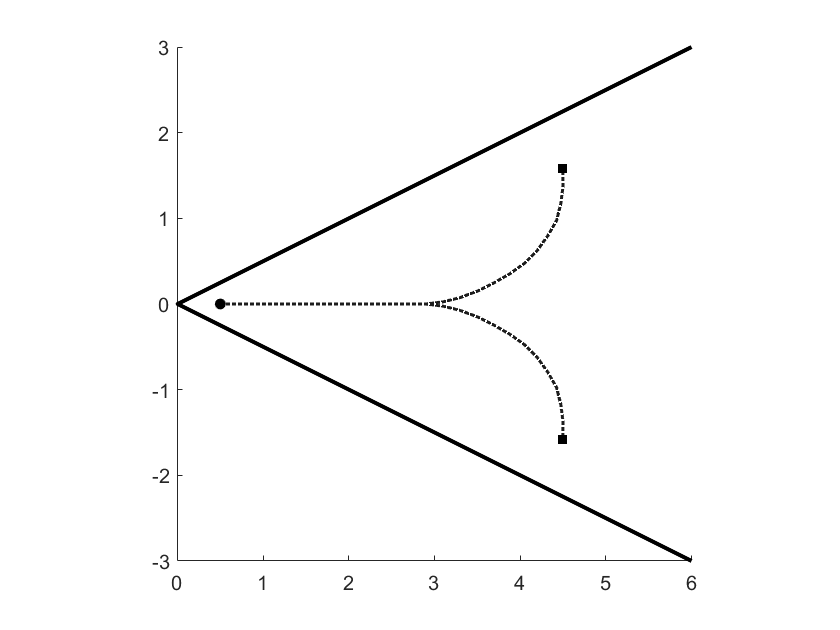}
    \caption{The case when $j=2$}
    \label{fig:convex_j2}
  \end{subfigure}
  \caption{Numerical parameters and corresponding bifurcation visualization for convex domain.}
  \label{fig:convex_analysis}
\end{figure}

Applying Algorithm \ref{alg:combined_algorithm} locally with setting $h=0.025$ and $m=8$ and applying Algorithm \ref{alg:bifurcation_point} to get $w_j$, we obtain the results given in Figure \ref{fig:convex_analysis}.

%\begin{table}[H]
%\begin{tabular}{@{}cccc@{}}
%\toprule
%$j$ & $w_j$ & $\epsilon_{\gamma_1}$ & $\epsilon_{\gamma_2}$ \\ \midrule
%0 & 2.9        & 0.1049        & 0.1049        \\
%1 & 2.9        & 0.0342        & 0.0342        \\
%2 & 2.9        & 0.0276        & 0.0276        \\
%3 & 2.9        & 0.0304        & 0.0304        \\
%4 & 2.9        & 0.033         & 0.033         \\
%5 & 3.1        & 0.02          & 0.02          \\
%6 & 3.1        & 0.0245        & 0.0245        \\
%7 & 3.1        & 0.0071        & 0.0071        \\
%8 & 3.1        & 0.0042        & 0.0042        \\
%9 & 3.1        & 2.2204e-16    & 2.2204e-16    \\ \bottomrule
%\end{tabular}
%\caption{The case $h=0.025$ and $m=8$.}
%\end{table}

%\begin{figure}[H]
%    \centering
%    \includegraphics[width=0.5\linewidth]{figs/bifurcation/convex domain/bif_h=0.025_m=8_j=2_grayscale.png}
%    \caption{The case when $j=2$}
%    \label{fig:placeholder}
%\end{figure}

Moreover, since we cannot define points on boundary, we use $z_0+10^{-8}i$ to replace $z_0$ for the case $j=0$. Based on the results, we conclude that the points $w_j$ are close to $y$, and $\epsilon_{\gamma_1}$ and $\epsilon_{\gamma_2}$ have an error less than $10^{-2}$ for $1 \leq j \leq 9$. 

\begin{figure}[H]
  \centering
  \begin{subfigure}{0.48\textwidth}
    \centering
    \includegraphics[width=\linewidth]{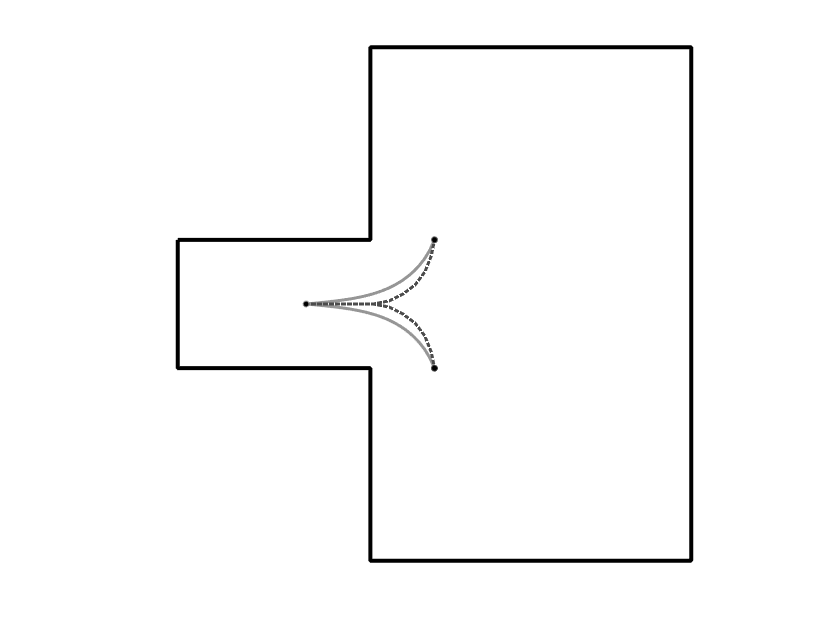}
    \caption{Connected curves are hyperbolic geodesics on $G$, and dotted curves are quasihyperbolic geodesics on $G$.}
    \label{fig: non_convex}
  \end{subfigure}
  \hfill
  \begin{subfigure}{0.48\textwidth}
    \centering
    \includegraphics[width=\linewidth]{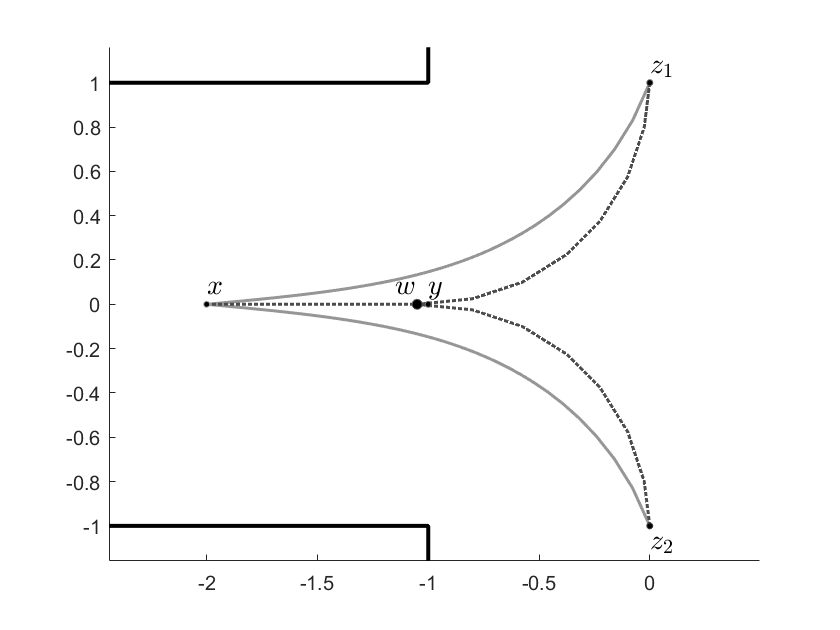}
    \caption{Let $x,y,z_1,z_2$ be points described in \ref{non-convex}, and curves are same as \ref{fig: non_convex} described.}
    \label{fig: non_convex_zoom_in}
  \end{subfigure}
  \caption{Bifurcation analysis in non-convex domain: comparison of hyperbolic and quasihyperbolic geodesics}
  \label{fig:non_convex_bifurcation}
\end{figure}

\begin{nonsec}{\bf Non-convex domains.} \label{non-convex}
    The following example is from \cite[Example 5.3]{AnJa}. Let $G$ be the polygonal domain bounded by the line segments connecting the points $(-4+i, -1+i, -1+4i, 4+4i, 4-4i, -1-4i,-1-i,-4-i)$, and let $x=-2,y=-1,z_1=i,z_2=-i$, where we set $y$ to be local bifurcation point and $w$ to be approximated bifurcation point. See Figure \ref{fig:non_convex_bifurcation}.
\end{nonsec}    

%\begin{figure}[H]
%  \centering
%  \includegraphics[width=0.5\linewidth]{figs/bifurcation/non-convex domain/non_convex_rst_grayscale.png}
%  \caption{Connected curves are hyperbolic geodesics on $G$, and dotted curves are quasihyperbolic geodesics on $G$.}
%  \label{fig: non_convex}
%\end{figure}

%\begin{figure}[H]
%  \centering
%  \includegraphics[width=0.5\linewidth]{figs/bifurcation/non-convex domain/h=0.025,m=8_grayscale.png}
%  \caption{$x,y,z_1,z_2$ are points described as above, and curves are same as \ref{fig: non_convex} described. This image applies DM with $h = 0.025, m=8$.}
%  \label{fig: non_convex_zoom_in}
%\end{figure}

In this experiment, we apply Algorithm \ref{alg:combined_algorithm} locally with several scales of refinement and apply Algorithm \ref{alg:bifurcation_point} to approximate $w$. Moreover, we define the error estimate as $\epsilon = |y-w|$ to describe the distance between $y$ and $w$. See Table \ref{table:error_terms_scales}.
\begin{table}[H]
\begin{tabular}{@{}ccccc@{}}
\toprule
$h$   & $m$ & $y$  & $w$      & $\epsilon$   \\ \midrule
0.05  & 2 & -1 & -0.8   & 0.2   \\
0.05  & 4 & -1 & -0.95  & 0.05  \\
0.025 & 2 & -1 & -0.8   & 0.2   \\
0.025 & 4 & -1 & -0.925 & 0.075 \\
0.025 & 8 & -1 & -1.05  & 0.05  \\
0.02  & 8 & -1 & -1.02  & 0.02  \\ 
0.01  & 2 & -1 & -0.76  & 0.24 \\ 
0.01  & 4 & -1 & -0.9  & 0.1 \\
0.01  & 8 & -1 & -0.96  & 0.04 \\
0.01  & 10 & -1 & -1  & 0 \\\bottomrule
\end{tabular}
\caption{Error estimates are dependent on different scales.}\label{table:error_terms_scales}
\end{table}

According to our experiments, when we refine the grid mesh by decreasing $h$ and increasing $m$, the approximated bifurcation point $w$ converges to $y$. In particular, it cannot be characterized by only considering the conformal density function.

%\textbf{Experiment: change $z_1,z_2$}
%The bifurcation point is dependent on choice %of$x,z_1,z_2$. 

\begin{comment}
    \section{Example 250930}
Let $G = R_3 \setminus L_3$, where $R_3$ and $L_3$ defined as \ref{sec: Lemma 5.1}. In this experiment, we will choose a set of complex values pairs $\{(z,w):z,w \in G\}$ with rules:
\begin{enumerate}
    \item $\Re(z) < \frac {\sqrt3} 2, \Re(w) > \frac {\sqrt3} 2$
    \item $|\Re(z)-\frac {\sqrt3} 2| = |\Re(w)-\frac {\sqrt3} 2|$
    \item $\Im(z) = \Im(w)$
\end{enumerate}
Therefore, we set $z = a+bi$ and $w = (\sqrt3-a)+bi$, with $-1< a < \frac {\sqrt 3} 2$ and $0 < b < 1$.

For each pair $(z,w)$, generate quasihyperbolic geodesic $\gamma_k(a,b)$ and find the behavior around $\frac {\sqrt3} 2$.

\begin{figure}[H]
    \centering
    \includegraphics[width=0.5\linewidth]{figs/250930/idea.png}
    \caption{idea for experiments}
    \label{fig:placeholder}
\end{figure}
\textbf{Table}

\textbf{Conclusion}
\end{comment}
%%%%%%%%%%%%%%%%%%%%%%%%%%%%
%%%%%%%%%%%%%%%%%%%%%%%%%%%%
%%%%%%%%%%%%%%%%%%%%%%%%%%%%

\section{Rectangular Domains}
%%%%%%%%%%%%%%%%%%%%%%%%%%%%
%%%%%%%%%%%%%%%%%%%%%%%%%%%%
%%%%%%%%%%%%%%%%%%%%%%%%%%%%

Consider the quasihyperbolic geodesics joining the point pairs $p_j$ and $q_j$, $j=1,2$, of different points. Then these two geodesics may have sub-arcs. This is closely related to the  medial axis of a domain \cite{hasto,Linden}, we study the following example:

Let $G = \{a+bi \; | \;a \in (-3,3),b \in (-1,1)\}$, and consider the points $z\in G$ with $\Re(z)=-2.9$, $\Im(z)=-0.8,-0.7, \ldots ,-0.1$, and $w = 2.5-0.5i$.

%\begin{figure}[H]
%    \centering
%    \includegraphics[width=0.5\linewidth]{figs/Rectangle Domain/domain.png}
%    \caption{Rectangular Domain}
%    \label{fig:placeholder}
%\end{figure}

\begin{figure}[H]
    \centering
    
    % 第一行
    \begin{subfigure}[b]{0.45\textwidth}
        \centering
        \includegraphics[width=\linewidth]{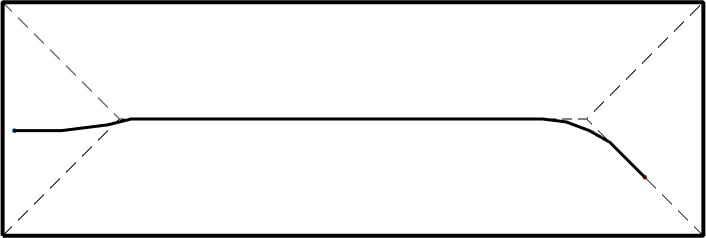}
        \caption{$\Im(z)=-0.1$}
        \label{fig:1a}
    \end{subfigure}
    \hfill
    \begin{subfigure}[b]{0.45\textwidth}
        \centering
        \includegraphics[width=\linewidth]{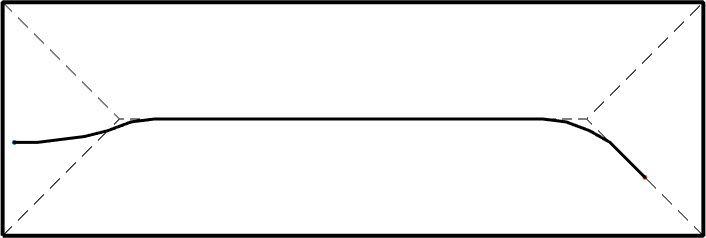}
        \caption{$\Im(z)=-0.2$}
        \label{fig:1b}
    \end{subfigure}
    
    % 第二行
    \begin{subfigure}[b]{0.45\textwidth}
        \centering
        \includegraphics[width=\linewidth]{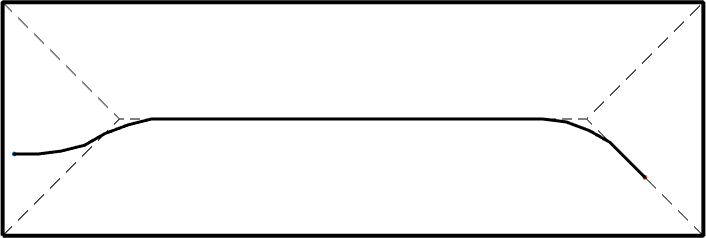}
        \caption{$\Im(z)=-0.3$}
        \label{fig:2a}
    \end{subfigure}
    \hfill
    \begin{subfigure}[b]{0.45\textwidth}
        \centering
        \includegraphics[width=\linewidth]{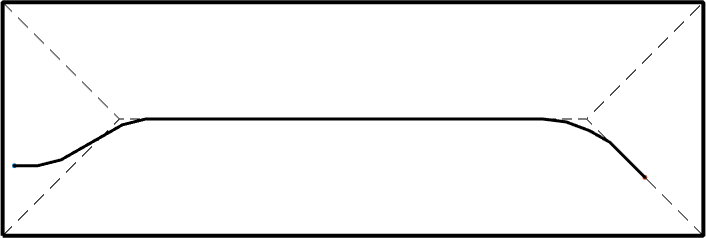}
        \caption{$\Im(z)=-0.4$}
        \label{fig:2b}
    \end{subfigure}
    
    % 第三行
    \begin{subfigure}[b]{0.45\textwidth}
        \centering
        \includegraphics[width=\linewidth]{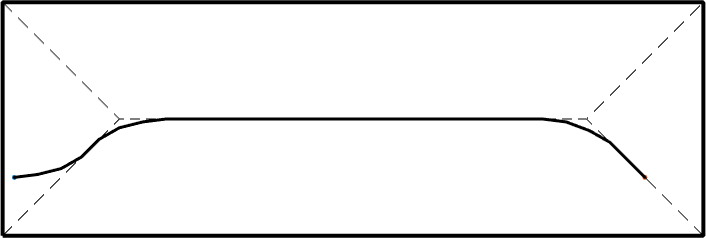}
        \caption{$\Im(z)=-0.5$}
        \label{fig:3a}
    \end{subfigure}
    \hfill
    \begin{subfigure}[b]{0.45\textwidth}
        \centering
        \includegraphics[width=\linewidth]{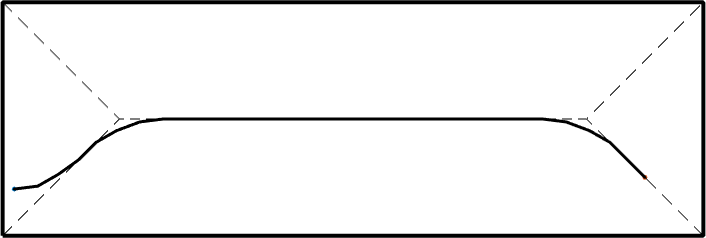}
        \caption{$\Im(z)=-0.6$}
        \label{fig:3b}
    \end{subfigure}
    
    % 第四行
    \begin{subfigure}[b]{0.45\textwidth}
        \centering
        \includegraphics[width=\linewidth]{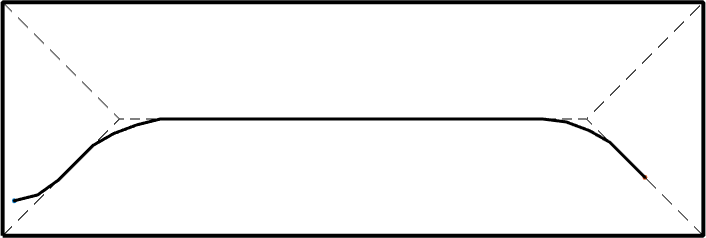}
        \caption{$\Im(z)=-0.7$}
        \label{fig:4a}
    \end{subfigure}
    \hfill
    \begin{subfigure}[b]{0.45\textwidth}
        \centering
        \includegraphics[width=\linewidth]{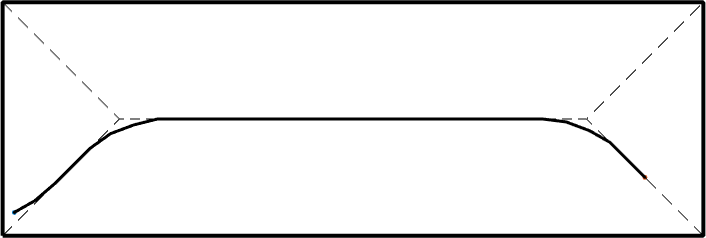}
        \caption{$\Im(z)=-0.8$}
        \label{fig:4b}
    \end{subfigure}
    
    \caption{Approximated quasihyperbolic geodesics in a rectangular domain.}
    \label{fig:collection_of_qh_j_in_G}
\end{figure}

Applying Algorithm \ref{alg:combined_algorithm} with $h=0.025, m=8$, we obtain the results in Figure \ref{fig:collection_of_qh_j_in_G}.

These results show  that the quasihyperbolic geodesic between $z$ and $w$ may have sub-arcs, which coincide with the medial axis. For some specific choices of $z$ in this example, we see that this sub-arc is a sub-segment of the angular bisector of the rectangle.

%%%%%%%%%%%%%%%%%%%%%%%%%%%%
%%%%%%%%%%%%%%%%%%%%%%%%%%%%
%%%%%%%%%%%%%%%%%%%%%%%%%%%%
\section{Asterisk-like Domains with $n$ Hands}
%%%%%%%%%%%%%%%%%%%%%%%%%%%%
%%%%%%%%%%%%%%%%%%%%%%%%%%%%
%%%%%%%%%%%%%%%%%%%%%%%%%%%%

\begin{figure}[H]
	\centering
	\begin{tikzpicture}[scale=1]
		%--BASIC SHAPE---------------------
		%------DEFINITIONS-----------------
		\def\nn{3}; % number of arms
		\def\lthree{1}; % fix l_3 = 1
		\def\sidelength{3}; % fix the length of arm's side
        \def\ltwo{2.5}; % fix the distance l_2: origin to p_j
		\def\bet{180/\nn}; % beta in degrees = (pi/n) * (180/pi) = 180/n
		\def\bdist{{\lthree/sin(\bet)}}; % distance of B_i from the origin
		\def\acdist{{(\sidelength+\lthree/tan(\bet))/cos(atan(\lthree/(\sidelength+\lthree/tan(\bet))))}}; % distance from origin to A_i and C_i
		\foreach \x in {1,...,\nn} \coordinate (a\x) at ({2*\x*\bet+atan(\lthree/(\sidelength+\lthree/tan(\bet)))}:\acdist); % corners A_i
		\foreach \x in {1,...,\nn} \coordinate (c\x) at ({2*\x*\bet-atan(\lthree/(\sidelength+\lthree/tan(\bet)))}:\acdist); % corners C_i
		\foreach \x in {1,...,\nn} \coordinate (b\x) at ({(2*\x-1)*\bet}:\bdist); % corners B_i
        \foreach \x in {1,...,\nn} \coordinate (P\x) at ({(2*\x)*\bet}:\ltwo); % points p_j
		\coordinate (o) at (0:0); % origin, needs to be defined like this for drawing angles
		\coordinate (p1) at (0:\bdist); % needed for drawing angles
		%-------LINES-----------------------
		\draw[thick] \foreach[evaluate=\x as \y using {mod(\x,\nn)+1}] \x in {1,...,\nn}{(b\x)--(c\x)--(a\x)--(b\y)}; % boundary
		\draw (o) circle(\bdist); % inner circle
		\draw (o)--(b1) node[midway, above] {$r$}; % radius of the inner circle with label
		%-------AXES------------------------
		\draw[->] ($(180:\acdist)+(180:1)$) -- ($(0:\acdist)+(0:1)$); % x-axis
		\draw[->] ($(270:\acdist)+(270:1)$) -- ($(90:\acdist)+(90:1)$); % y-axis
		%-------NODES-----------------------
		\draw \foreach \x in {1,...,\nn} \foreach \y in {a,b,c,P} {
			(\y\x) node[circle, fill=black, inner sep=0pt, minimum width=5pt] {} }; % corner nodes A_i, B_i, C_i and points p_j
		\draw (o) node[circle, fill=black, inner sep=0pt, minimum width=5pt] {}; % origin
		%--END OF BASIC SHAPE---------------
		%-----------------------------------
		%--CORNER LABELS--------------------
		% For large values of \nn, the coefficients for the coordinates need to be adjusted. These coefficients work well for \nn = 3,...,5.
		\draw \foreach[evaluate=\x as \y using {int(mod(\x,\nn)+1)}] \x in {1,...,\nn}{($1.12*(a\x)$) node[circle, fill=white, inner sep= 0pt] {$A_{\y}$}}; % labels of A_i: what is shown as A_2 is in fact (a2) and A_1 is (a\nn)
		\draw \foreach[evaluate=\x as \y using {int(mod(\x,\nn)+1)}] \x in {1,...,\nn}{($1.12*(c\x)$) node[circle, fill=white, inner sep= 0pt] {$C_{\y}$}}; % labels of C_i
        \draw \foreach[evaluate=\x as \y using {int(mod(\x,\nn)+1)}] \x in {1,...,\nn}{($1.15*(P\x)$) node[circle, fill=white, inner sep= 0pt] {$p_{\y}$}}; % labels of p_j
		\draw \foreach \x in {2,...,\nn}{($1.45*(b\x)$) node[circle, fill=white, inner sep= 0pt] {$B_{\x}$}}; % labels of B_i, i=2,...,n
        \draw ($1.5*(b1)-(0.2,0)$) node[circle, fill=white, inner sep= 0pt] {$B_{1}$}; % label of B_1
		\draw (o) node[above left] {$O$}; % origin
		%--LENGTH LABELS AND BRACES---------
		\draw [thick, decorate, decoration = {calligraphic brace,raise=3pt,amplitude=5pt}] (0:{\sidelength+\lthree/tan(\bet)}) -- (a\nn) node[pos=0.5,left=6pt,black] {$\ell_3$}; % l_3
		\draw [thick, decorate, decoration = {calligraphic brace,raise=3pt,amplitude=5pt}] (b\nn)--(c\nn) node[pos=0.5,above=6pt,black] {$\ell_1$}; % l_1
		%--ANGLES---------------------------
		\draw pic [draw, angle radius=9 pt] {angle=p1--o--b1}; % beta
        \draw ($(o)+(.5,.25)$) node[] {$\beta$}; % beta label
		\draw pic [draw, angle radius=8 pt] {angle=a\nn--b1--c1}; % theta
        \draw ($(b1)+(.5,.25)$) node[] {$\theta$}; % theta label
		\draw pic [draw, angle radius=7 pt] {right angle=b1--a\nn--c\nn}; % right angle A_1
		\draw pic [draw, angle radius=7 pt] {right angle=a\nn--c\nn--b\nn}; % right anlge C_1
	\end{tikzpicture}
    \caption{A 3-asterisk. }\label{fig:astarisk}
\end{figure}
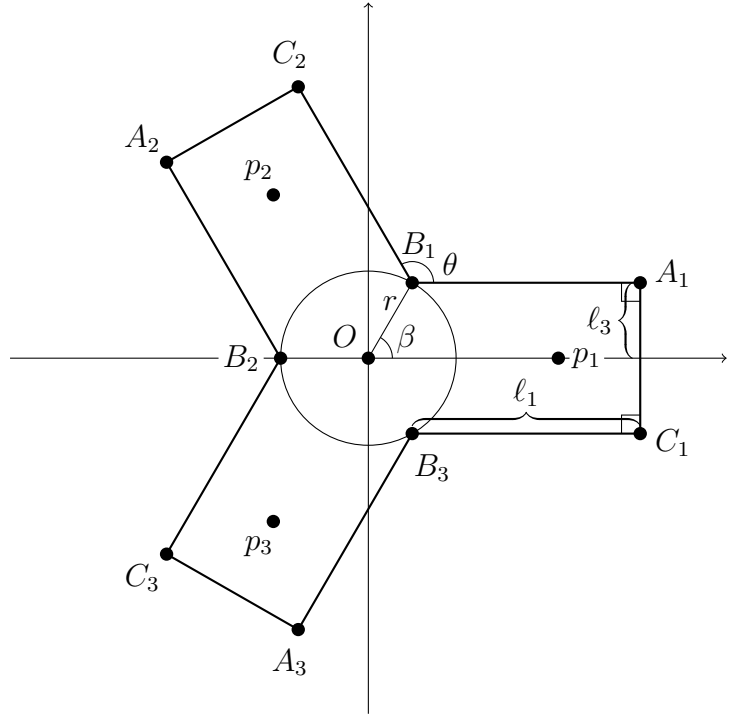

Choose initial values $\ell_1,\ell_2,\ell_3>0$, and $n \in \mathbb{N}$ such that $n \geq 3$. Let $\theta = \frac {2\pi} n$, $\alpha = \frac \pi 2- \frac \pi n$, $\beta = \frac \pi n$, $r = \frac {l_3} {\sin{\beta}}$, $\rho = \sqrt{\ell_1^2+r^2+2\ell_1r\cos{\beta}}$, $\phi = \arctan{\frac {\ell_3} {r\cos{\beta}+\ell_1}}$. We call the polygon $G$ with $n\ge3$ counterclockwise vertices given by
\[
\{C_1,A_1,B_1,C_2,\dots,C_j,A_j,B_j,\dots, B_{n-1},C_n,A_n,B_n:1 \leq j \leq n\}
\]
with $A_j = \rho \cdot e^{i(\phi+j\theta)}, B_j = r\cdot e^{i(\beta+j\theta)}, C_j = \rho \cdot e^{i(-\phi+j\theta)}$, an $n$-Asterisk. The case $n=3$ is illustrated in Figure \ref{fig:astarisk}.

Let $(p_1,p_2,\dots,p_n,p_1)$ be points to be connected with hyperbolic geodesics and quasihyperbolic geodesics, where $p_j = (\Re(A_1)-\ell_2) \cdot e^{j \theta}$.

Our aim is to find the inscribed circle of hyperbolic polygon $D_\rho = \{\gamma_\rho(p_j,p_{j+1})\}$ and quasihyperbolic polygon $D_k = \{\gamma_k(p_j,p_{j+1})\}$, comparing their radius, and to observe their geometry. In our experiment, the radii of their inscribed circles are calculated from the formulas:
\[ 
\left\{\begin{array}{rcll}
r_\rho &=& \min_{1\leq j \leq n} \{ \min \{|z_k|:z_k \in \gamma_j\} \},& \gamma_j = \gamma_\rho(p_j,p_{j+1}),\\
r_k &=& \min_{1\leq j \leq n} \{ \min \{|z_k|:z_k \in \gamma_j\} \}, & \gamma_j = \gamma_k(p_j,p_{j+1}).
\end{array}\right.
\]

\begin{figure}[H]
  \centering
  \begin{subfigure}{0.40\textwidth}
    \centering
    \includegraphics[width=\linewidth]{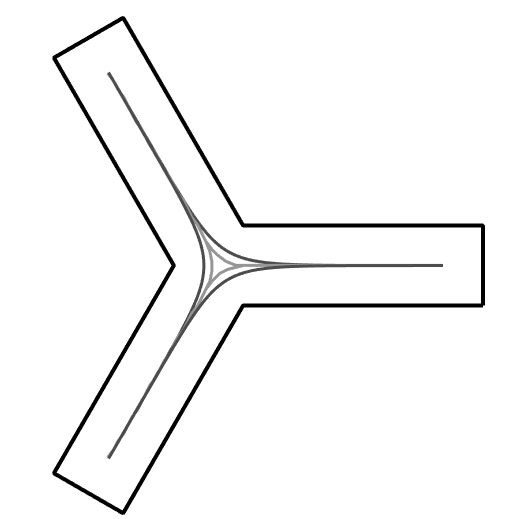}
    \caption{The dark curve is a hyperbolic geodesic, and the light curve is a quasihyperbolic geodesic.}
    \label{fig: n=3_grayscale}
  \end{subfigure}
  \hfill
  \begin{subfigure}{0.58\textwidth}
    \centering
    \includegraphics[width=\linewidth]{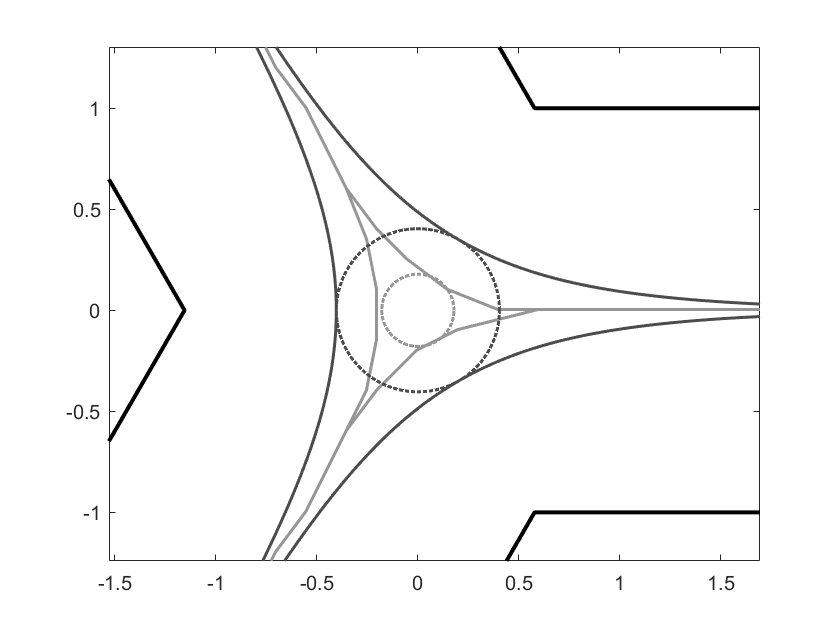}
    \caption{The dark dotted circle is the largest inscribed circle of a hyperbolic polygon $D_\rho$, and the light dotted circle is the largest inscribed circle of quasihyperbolic polygon $D_k$.}
    \label{fig: n=3_grayscale_zoom_in}
  \end{subfigure}
  \caption{Comparison of hyperbolic and quasihyperbolic geodesics in the $n$-asterisk for $n=3$.}
  \label{fig:n3_geodesics_comparison}
\end{figure}

In the experiment, we initially set $\ell_1 = 6$, $\ell_2=\ell_3=1$ and apply Algorithm \ref{alg:combined_algorithm} with the parameter values $h=0.05$ and $ m=4$. The case $n=3$ is illustrated in Figure \ref{fig:n3_geodesics_comparison}.

%\begin{figure}[H]
%  \centering
%  \includegraphics[width=0.5\linewidth]{figs/n-astarisk/n=3_grayscale.png}
%  \caption{The dark curve is hyperbolic geodesic, and the light curve is quasihyperbolic geodesic.}
%  \label{fig: n=3_grayscale}
%\end{figure}

%\begin{figure}[H]
%  \centering
%  \includegraphics[width=\linewidth]{figs/n-astarisk/n=3,with_circle_grayscale.png}
%  \caption{The dark dotted circle is inscribed circle of $D_\rho$, and the light dotted circle is inscribed circle of $D_k$,}
%  \label{fig: n=3_grayscale_zoom_in}
%\end{figure}

We obtain the following results by changing $n$ from $3$ to $9$.

%\textbf{change numbers to $10^4$}
%\begin{table}[H]
%\begin{tabular}{@{}cccccc@{}}
%\toprule
%$n$ & $Re(p_1)$ & $r$      & $r_\rho$ & $r_k$   & $\frac {r_\rho} {r_k}$ \\ \midrule
%3 & 5.5774       & 1.1547 & 0.40365  & 0.17860 & 2.260145182  \\
%4 & 6.0000       & 1.4142 & 0.73942  & 0.42426 & 1.742838932  \\
%5 & 6.3764       & 1.7013 & 1.0604   & 0.70496 & 1.50419923   \\
%6 & 6.7321       & 2.0000 & 1.3772   & 0.99975 & 1.377591899  \\
%7 & 7.0765       & 2.3048 & 1.6928   & 1.3060  & 1.296227781  \\
%8 & 7.4142       & 2.6131 & 2.0081   & 1.6169  & 1.241989911  \\
%9 & 7.7475       & 2.9238 & 2.3235   & 1.9082  & 1.217612586  \\ \bottomrule

%\end{tabular}
%\caption{$l_1 = 6$, $l_2=l_3=1$, with $h=0.05, m=4$}
%\end{table}

\begin{table}[H]
\begin{tabular}{@{}cccccc@{}}
\toprule
$n$ & $\Re(p_1)$ & $r$ & $r_\rho$ & $r_k$ & $\frac {r_\rho} {r_k}$ \\ \midrule
3 & 5.5774 & 1.1547 & 0.40365 & 0.17860 & 2.2601 \\
4 & 6.0000 & 1.4142 & 0.73942 & 0.42426 & 1.7428 \\
5 & 6.3764 & 1.7013 & 1.0604 & 0.70496 & 1.5042 \\
6 & 6.7321 & 2.0000 & 1.3772 & 0.99975 & 1.3776 \\
7 & 7.0765 & 2.3048 & 1.6928 & 1.3060 & 1.2962 \\
8 & 7.4142 & 2.6131 & 2.0081 & 1.6169 & 1.2420 \\
9 & 7.7475 & 2.9238 & 2.3235 & 1.9082 & 1.2176 \\ \bottomrule
\end{tabular}
\caption{Experimental results for $\ell_1 = 6$, $\ell_2=\ell_3=1$, with $h=0.05, m=4$.}
\end{table}

%\textbf{subfigures n=4,5,6,7,8,9}
\begin{figure}[htbp]
    \centering
    
    \begin{subfigure}[b]{0.45\textwidth}
        \centering
        \includegraphics[width=\linewidth]{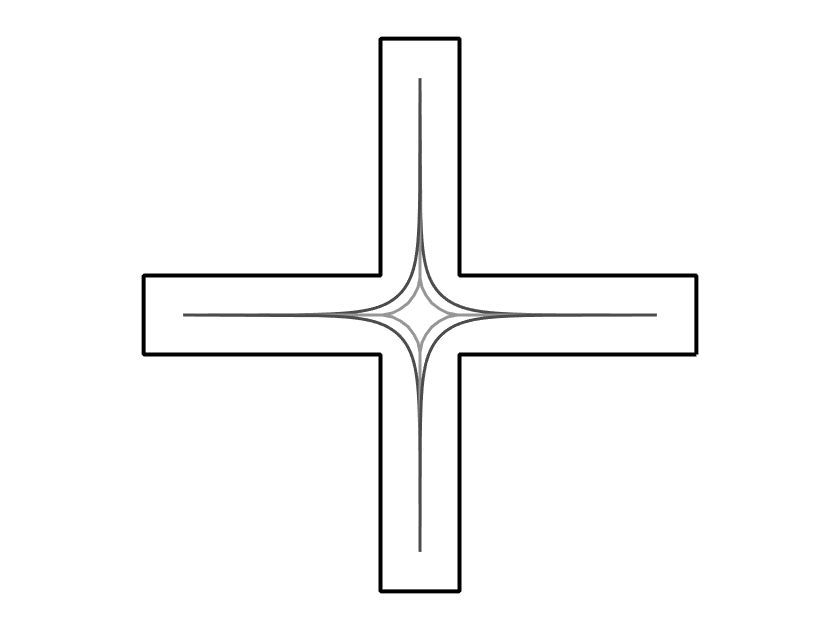}
        \caption{$n=4$}
        \label{fig:4-astarisk}
    \end{subfigure}
    \hfill
    \begin{subfigure}[b]{0.45\textwidth}
        \centering
        \includegraphics[width=\linewidth]{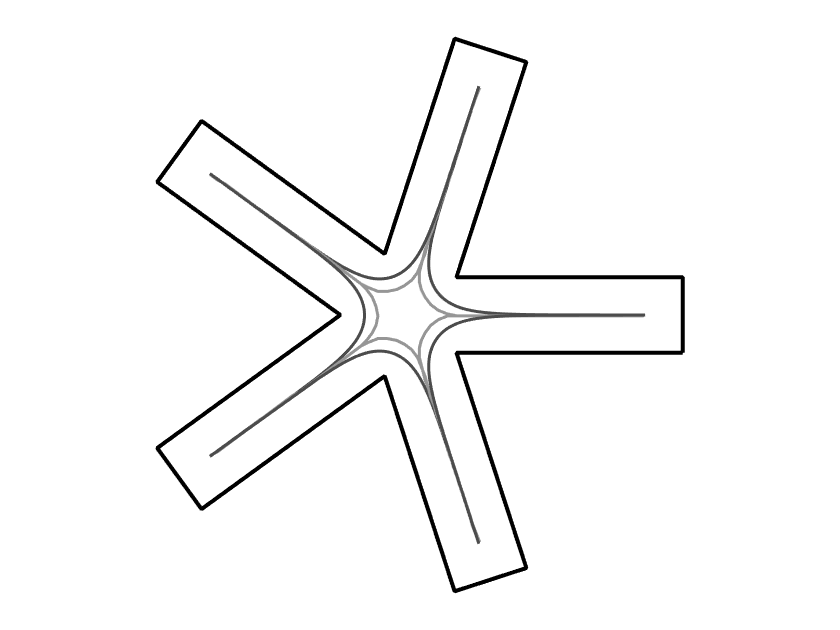}
        \caption{$n=5$}
        \label{fig:5-astarisk}
    \end{subfigure}
    
    \begin{subfigure}[b]{0.45\textwidth}
        \centering
        \includegraphics[width=\linewidth]{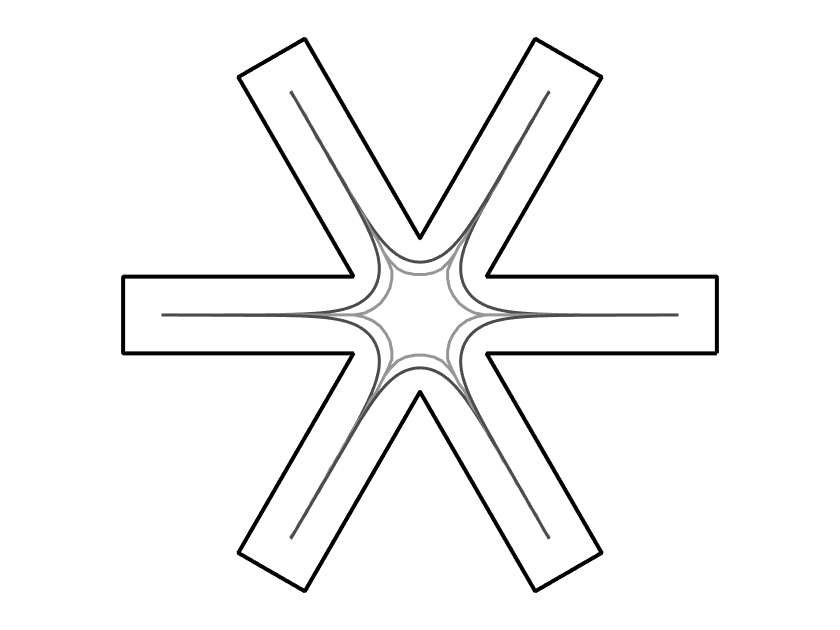}
        \caption{$n=6$}
        \label{fig:6-astarisk}
    \end{subfigure}
    \hfill
    \begin{subfigure}[b]{0.45\textwidth}
        \centering
        \includegraphics[width=\linewidth]{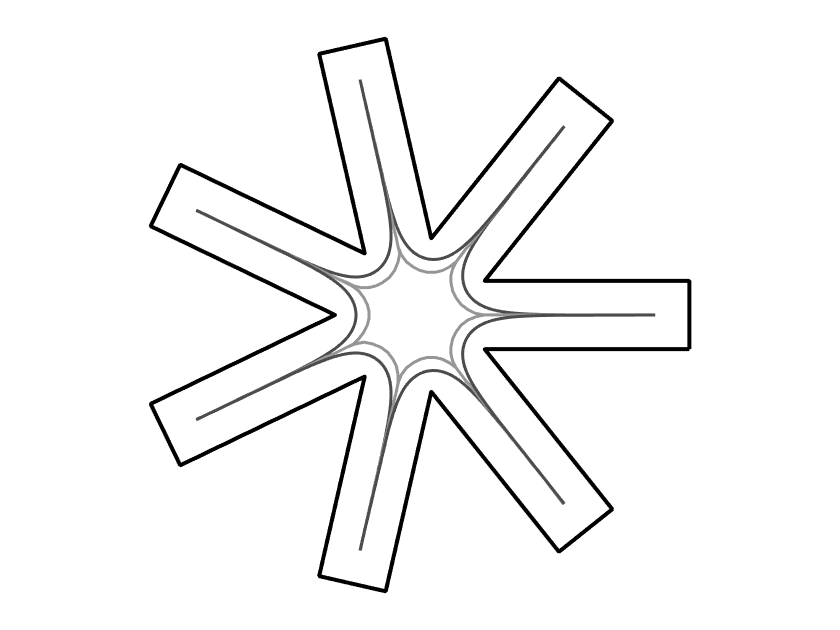}
        \caption{$n=7$}
        \label{fig:7-astarisk}
    \end{subfigure}
    
    \begin{subfigure}[b]{0.45\textwidth}
        \centering
        \includegraphics[width=\linewidth]{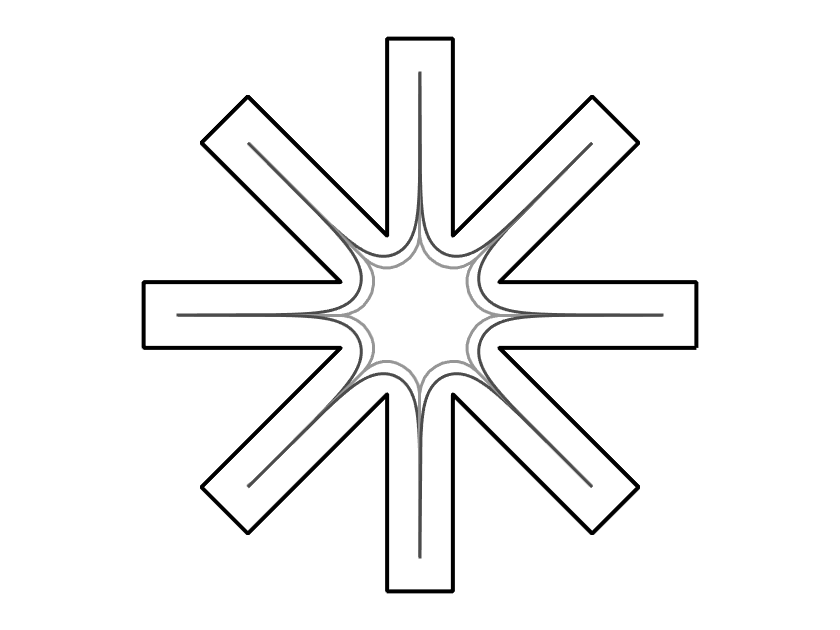}
        \caption{$n=8$}
        \label{fig:8-astarisk}
    \end{subfigure}
    \hfill
    \begin{subfigure}[b]{0.45\textwidth}
        \centering
        \includegraphics[width=\linewidth]{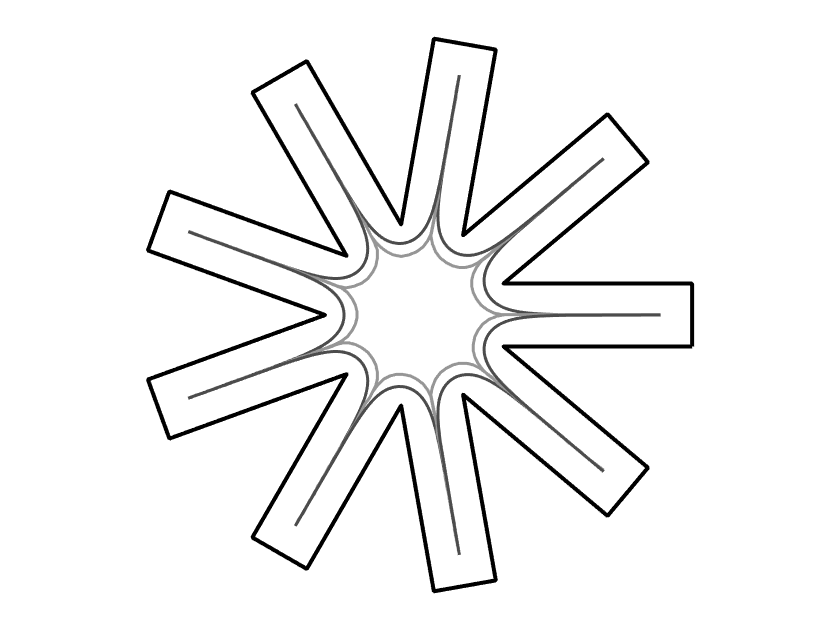}
        \caption{$n=9$}
        \label{fig:9-astarisk}
    \end{subfigure}

    \caption{Hyperbolic and quasihyperbolic polygons in $n$-asterisks.}
    \label{fig:collection_qh_n-astarisks}
\end{figure}

Based on the above results, we conclude that the inscribed circle of a hyperbolic polygon must contain the inscribed circle of quasihyperbolic polygon. Approximations of geodesics are illustrated in Figure \ref{fig:collection_qh_n-astarisks}.

\section{Remarks on Dijkstra's and A* Algorithms}
%%%%%%%%%%%%%%%%%%%%%%%%%%%%
%%%%%%%%%%%%%%%%%%%%%%%%%%%%
%%%%%%%%%%%%%%%%%%%%%%%%%%%%
Dijkstra's algorithm is a so-called {\it uniform cost search} algorithm, and it traverses every node of the given graph, calculates lengths of  all possible paths, which makes it a low efficiency process. With respect of this question, F.C. Chrislock \cite{fc} came up with the idea of applying a {\it heuristic algorithm} such as A* to find shortest paths, i.e., geodesics.

The difficulty is to find the heuristic function $h(n)$, which is done just like what we do with Dijkstra's method, by applying hyperbolic or quasihyperbolic metric as our weight function. However, since we do not know the explicit formula for the quasihyperbolic metric, it is not possible to construct such a mesh with quasihyperbolically uniform  discretization step $h$. It seems that Chrislock was also aware of this phenomenon, so he avoided searching for a formula and mapped the conflict region conformally onto to a disk instead of using the Euclidean distance as heuristic function.

\begin{figure}[htbp]
    \centering
    
    \begin{subfigure}[b]{0.45\textwidth}
        \centering
        \includegraphics[width=\linewidth]{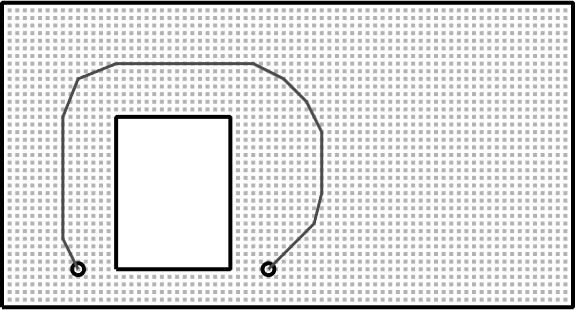}
        \caption{$h_1 = 0.2$}
        \label{fig: large curve}
    \end{subfigure}
    \hfill
    \begin{subfigure}[b]{0.45\textwidth}
        \centering
        \includegraphics[width=\linewidth]{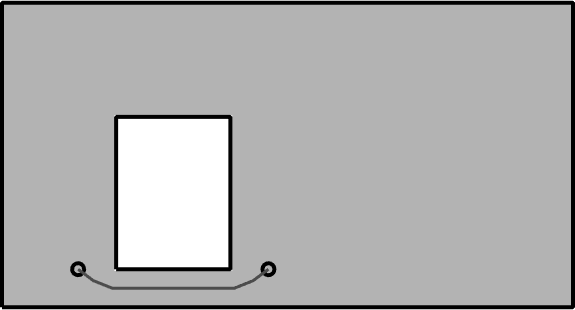}
        \caption{$h_2 = 0.1$}
        \label{fig: short curve}
    \end{subfigure}
    
    \caption{Approximated quasihyperbolic geodesics with different $h$.}
    \label{fig:qh_j-doubly-connected}
\end{figure}

One possible approach is to apply Algorithm \ref{alg:combined_algorithm} to obtain a curve first, then apply A* around a neighborhood of the curve to refine the results. This reduces
the size of the graph significantly
and therefore accelerates the computation. However, when designing the experiments, we considered the following case:

Let $G = R \setminus P$, where $R = \{a + bi \; | \; a \in [0,15], b \in [0,8]\}$ and $P = \{a + bi \; | \; a \in [3,6], b \in [1,5]\}$, where $p_1 = 2 + i$ and $p_2 = 7 + i$ are the points to be connected, see Figure \ref{fig:qh_j-doubly-connected}. We apply Dijkstra's method on $G$ with different values of $h$ and $m$ to obtain the approximation of a quasihyperbolic geodesic between $p_1$ and $p_2$. When setting $m=4$, the shapes of the geodesics are totally different with $h_1=0.2$ and $h_2=0.1$, respectively.

Moreover, $k_{G_{h_1,m}}^*(p_1,p_2) \approx 11.22949$ and $k_{G_{h_2,m}}^*(p_1,p_2) \approx 10.04478$. This example shows that refinement of the mesh may change the route taken by the approximated geodesic.  Therefore, one needs to apply Dijkstra's method on $G$ globally in order to handle the topology of multiply connected domains.
%%%%%%%%%%%%%%%%%%%%%%
%%%%%%%%%%%%%%%%%%%%%%
%%%%%%%%%%%%%%%%%%%%%%
\section{Conclusions}

We have demonstrated that applying Algorithm \ref{alg:combined_algorithm} can be used to approximate quasihyperbolic geodesics on simply and multiply connected domains, including domains with punctures. In these cases we are able to obtain numerical results that are compatible with earlier theoretical results.

Our methods can be used to experimental study of many problems in geometric function theory, including theoretical properties of quasihyperbolic geodesics and balls, and triangles that play an important role in the question of Gromov hyperbolicity of these metrics. Moreover, these methods can be applied not only for hyperbolic and quasihyperbolic metrics, but for other weighted metrics as well, including Riemannian metrics on surfaces. Therefore  we anticipate many potential applications.

%%%%%%%%%%%%%%%%%%%%%%
%%%%%%%%%%%%%%%%%%%%%%
%%%%%%%%%%%%%%%%%%%%%%
%%%%%%%%%%%%%%%%%%%%%%%%%%%%%%%%%%%%%%%%%%%%%%%%%%%%%%%
\subsection*{Acknowledgments}
S.G. and A.R. were supported by Natural Science Foundation of Guangdong Province (No. 2024A1515010467), Li Ka~Shing Foundation GTIIT-STU Joint Research Grant (No. 2024LKSFG06), and GTIIT Education Foundation.

\medskip

%%%%%%%%%%%%%%%%%%%%%%%%%%%%%%%%%%%%%%%%%%%%%%%%%%%%%%%
\bibliographystyle{siam}

\bibliography{references}

\end{document}